\documentclass[a4paper,12pt]{article}
\usepackage{graphicx} 

\usepackage{amsmath}
\usepackage{multirow,bigdelim}
\usepackage{amscd}
\usepackage{amssymb,amsfonts,amsmath,amsthm}
\usepackage{verbatim}

\usepackage{amsmath}
\usepackage{amsthm}
\usepackage{amssymb}

\usepackage{latexsym}
\usepackage{array}
\usepackage{enumerate}
\usepackage[initials,alphabetic]{amsrefs}
   \BibSpec{arXiv}{
   +{}{\PrintAuthors}{author}
   +{,}{ \textit}{title}
   +{,}{ }{date}
   +{,}{ arXiv:}{eprint}}

\numberwithin{equation}{section}

\usepackage{amsmath,amssymb,mathrsfs,amsthm}
\usepackage[all]{xy}
\usepackage{graphicx}
\usepackage[normalem]{ulem}
\usepackage{ascmac}
\usepackage{mathtools}

\newcommand{\rhom}{\mathrm{R}{\mathcal{H}}om}

\newcommand{\Mod}{\mathrm{Mod}}

\newcommand{\Perv}{\mathrm{Perv}}
\newcommand{\shom}{{\mathcal{H}}om}

\newcommand{\rsect}{{\mathrm{R}}\Gamma}
\newcommand{\gr}{\mathrm{gr}}
\newcommand{\SD}{\mathcal{D}}
\newcommand{\SE}{\mathcal{E}}

\newcommand{\SM}{\mathcal{M}}
\newcommand{\SN}{\mathcal{N}}
\newcommand{\SO}{\mathcal{O}}

\newcommand{\calS}{\mathcal{S}}

\newcommand{\SC}{\mathcal{C}}

\newcommand{\SI}{\mathcal{I}}

\newcommand{\bbC}{\mathbb{C}}
\newcommand{\bbZ}{\mathbb{Z}}
\newcommand{\bbR}{\mathbb{R}}
\newcommand{\bbQ}{\mathbb{Q}}

\newcommand{\rmR}{\mathrm{R}}
\newcommand{\rmt}{\mathrm{t}}
\newcommand{\rmI}{\mathrm{I}}

\newcommand{\bfD}{\mathbf{D}}

\newcommand{\Sol}{\mathrm{Sol}}
\newcommand{\DR}{\mathrm{DR}}
\newcommand{\Db}{\mathbf{D}^\mathrm{b}}
\newcommand{\Dbc}{{\bf D}_{\mathrm{c}}^{\mathrm{b}}}

\newcommand{\tl}[1]{\widetilde{#1}}

\newcommand{\simto}{\overset{\sim}{\longrightarrow}}

\newcommand{\Modh}{\mathrm{Mod}_{\mathrm{h}}}
\newcommand{\Modcoh}{\mathrm{Mod}_{\mbox{\scriptsize coh}}}

\newcommand{\Modrh}{\mathrm{Mod}_{\mbox{\scriptsize rh}}}

\newcommand{\Dbhol}{{\mathbf{D}}^{\mathrm{b}}_{\mathrm{h}}}
\newcommand{\Dbrh}{{\mathbf{D}}^{\mathrm{b}}_{\mbox{\scriptsize rh}}}

\newcommand{\Dotimes}{\overset{D}{\otimes}}

\DeclareMathOperator{\ch}{ch}
\DeclareMathOperator{\CCyc}{CC}
\DeclareMathOperator{\supp}{supp}
\DeclareMathOperator{\RH}{RH}

\DeclareMathOperator{\Eu}{Eu}
\DeclareRobustCommand{\longtwoheadrightarrow}{\relbar\joinrel\twoheadrightarrow} 
\newcommand{\op}{\mathrm{op}}
\newcommand{\reg}{\mathrm{reg}}
\newcommand{\codim}{\operatorname{codim}}
\usepackage{calligra}
\newcommand{\orsh}{\text{\raisebox{0.2mm}{\Large\textcalligra{or}}}\!\;} 
\newcommand{\Sa}{\mathcal{A}}
\newcommand{\Dbrc}{\bfD_{\bbR\text{-}\mathrm{c}}^{\mathrm{b}}}
\renewcommand{\Re}{\operatorname{Re}}

\DeclareMathOperator{\MS}{SS}
\newcommand{\Eb}{\mathbf{E}^\mathrm{b}}

\DeclareRobustCommand{\longhookrightarrow}{\lhook\joinrel\longrightarrow}

\newtheorem{theorem}{Theorem}[section]
\newtheorem{conjecture}[theorem]{Conjecture}

\theoremstyle{definition}
\newtheorem{definition}[theorem]{Definition}
\theoremstyle{remark}

\newtheorem{example}[theorem]{\sc Example}

\theoremstyle{plain}
\newtheorem{theorem and definition}[theorem]{Theorem and Definition}

\title{The Legend of Masaki Kashiwara \\
and Algebraic Analysis\\[0.5em]
\large To the memory of Professor Mikio Sato}
\author{Kiyoshi TAKEUCHI 
\footnote{Mathematical Institute, Tohoku University, 
Aramaki Aza-Aoba 6-3, Aobaku, Sendai, 980-8578, Japan. 
E-mail: takemicro@nifty.com} }
\date{\today}

\begin{document}

\maketitle

\begin{abstract}
This survey paper offers a concise introduction to Kashiwara's work 
on $\SD$-modules, microlocal analysis and related subjects. In this way, 
we explain his role in the development of algebraic analysis. 
\end{abstract}

\section{Introduction}

As one of many mathematicians who have benefited greatly from the works and 
ideas of Masaki Kashiwara, I am very pleased to 
take on the role of conveying my enthusiasm for 
his remarkable achievements in mathematics to the reader. 
Because of the rapid development of different fields of mathematics, 
nowadays it is becoming harder and harder to move 
across several subjects, as historical giants such as 
Euler, Gauss, and Hilbert did centuries ago. 
Nevertheless, this is not the case with Masaki Kashiwara. 

One of the keys to his success lies in his 
systematic use of category theory. 
Arguments involving categories and functors 
may appear too general at first glance, but they allow us to transmit ideas 
between fields in a highly efficient way. 
In my opinion, Kashiwara is one of the very few mathematicians 
who have succeeded in making full use of this strength of category theory 
to obtain decisive results in various distinct branches of 
mathematics. 

In particular, he developed 
the theory of derived categories, introduced by 
Grothendieck and Mikio Sato in the 1960's, to 
its full potential, so that it can be applied to a broad 
range of problems. As a consequence, in the monumental 
paper \cite{SKK}, Sato, Kawai and Kashiwara constructed 
a general theory for systems of linear PDEs, of a kind that had never  
been envisioned before.  
Moreover in \cite{Kas6} and \cite{Kas9}, 
Kashiwara established the Riemann-Hilbert correspondence, 
by which problems in analysis (regular holonomic 
$\SD$-modules) can be translated into problems 
in geometry (perverse sheaves), and used it 
to  settle the Kazhdan-Lusztig conjecture in spectacular fashion 
with Brylinski in \cite{BK}. 

In this way, Kashiwara and  his coauthors 
tied together analysis and algebra, 
which had seemed  far removed from each other, 
and created a completely new field called ``algebraic 
analysis". Soon after its emergence in the 1970's,  
it became clear that the methods and ideas developed in this field 
were extremely powerful and applicable to many areas in mathematics, 
such as algebraic geometry, representation theory, 
singularity theory, number theory, and mathematical 
physics. Just as analysis was profoundly 
influenced by functional analysis a century ago, 
algebraic analysis transformed it once again at a fundamental level. 

In this revolution of mathematics, 
Kashiwara and others also used the theory of derived 
categories to create new fields and problems, 
greatly enlarging our vision of mathematics. 
In particular, Kashiwara's books \cite{KS2} and \cite{KS3} 
on derived categories, written in collaboration 
with Schapira, have had an enormous impact on mathematics 
and are widely regarded as foundational texts. 
Indeed, after Kashiwara, the mathematical landscape 
has been profoundly changed. Like great art or music, 
everybody who gets in touch with Kashiwara's work 
is struck by its beauty and elegance. In this sense, he has enriched human 
culture itself, 
gifting us a remarkable collection of intellectual masterpieces. 
Undoubtedly, his influence will continue to spread over the generations. 

The aim of this paper is to give a concise introduction to 
this rapidly growing area of mathematics, with 
special emphasis on Kashiwara's works, and assuming only that the 
reader is familiar with some basic definitions from sheaf theory. 
Given the extraordinary breadth of Kashiwara's contributions, 
I have selected only a few important topics that I am 
most familiar with. Each of these will be discussed separately in the 
subsequent sections. Other important topics will be treated by the other 
authors in this special volume dedicated to 
Masaki Kashiwara.

\bigskip
\noindent{\bf Acknowledgement:} 
The author thanks Ren Fernandes, Kazuki Kudomi and 
Hajime Tsuchiya for reading the manuscript 
carefully. He is also grateful to Professor 
Alan Stapledon and the anonymous referee 
for many useful suggestions on this paper.

\section{$\SD$-modules}
The theory of $\SD$-modules was initiated by Kashiwara \cite{Kas12} 
and Bernstein \cite{Ber1} independently.
Whereas Kashiwara's master thesis \cite{Kas12} treats analytic $\SD$-modules 
on complex manifolds, in Bernstein's paper \cite{Ber1} 
algebraic $\SD$-modules i.e. modules over Weyl algebras are treated. 
Here we mainly restrict ourselves to Kashiwara's theory of 
analytic $\SD$-modules. 
For the details of this section, see e.g. 
\cite{Bj1}, \cite{Bj}, \cite{ALD}, \cite{H-T-T}, 
\cite{K-book}, \cite{KS16}, \cite{Sc1}. 
Let $X$ be a complex manifold and denote by $\SO_X$ the sheaf 
of holomorphic functions on it. 
Then we define a sheaf $\SD_X$ of (non-commutative) rings on $X$ 
by assigning to each open subset $U \subset X$ 
the ring of differential operators on $U$ with coefficients in $\SO_X$. 
More precisely, we construct $\SD_X$ as a subsheaf of rings of 
$\SE nd_{\bbC_X} ( \SO_X )= \shom(\SO_X,\SO_X)$ and for two sections 
$P,Q \in \SD_X(U)$ on an open subset $U \subset X$ 
their product $PQ \in \SD_X(U)$ is defined to be the composition $P \circ Q$ 
of the $\bbC_X$-linear sheaf homomorphisms 
$\SO_U \stackrel{Q\times}\longrightarrow \SO_U \stackrel{P\times}\longrightarrow \SO_U$. 
If the open subset $U \subset X$ is non-emptry, connected and 
endowed with a holomorphic coordinate 
$x=(x_1,\ldots,x_n)$, then we have 
\begin{equation*}
\SD_X(U) = \biggl\{\sum_{\alpha \in \bbZ_+^n} a_\alpha(x)
\partial_x^\alpha \,\biggm\vert\, a_\alpha(x) \in \SO_X(U) \biggr\},
\end{equation*}
where we set $\bbZ_+ \coloneqq \{m\in \bbZ \,\vert\, m\geq0 \}$, 
the sums $\sum_{\alpha \in \bbZ_+^n}$ are finite, and for 
$\alpha=(\alpha_1,\ldots,\alpha_n) \in \bbZ_+^n$ the standard 
notation $\partial_x^\alpha = \partial_{x_1}^{\alpha_1} 
\cdots \partial_{x_1}^{\alpha_n}$ is used. 
In this situation, the products in the ring $\SD_X(U)$ are described 
by the classical Leibniz rule. 
For example, we have the relation 
$[\partial_{x_i},x_j] = \delta_{ij} \ (1\leq i,j \leq n)$. 
For a differential operator $P = \sum_{\alpha \in \bbZ_+^n} a_\alpha(x)
\partial_x^\alpha \in \SD_X(U)$ on such $U \subset X$ 
we set
\begin{equation*} 
\textrm{ord}P \coloneqq \textrm{max} \{\vert\alpha\vert = \alpha_1+
\cdots+\alpha_n \,\big\vert\, a_\alpha(x) \not\equiv 0\}
\end{equation*}
and call it the order of $P$. 
It turns out that this local definition does not depend on the choice of the local coordinate 
$x = (x_1,\ldots,x_n)$ and hence by setting
\begin{equation*}
F_i\SD_X \coloneqq \{P \in \SD_X \, \bigm\vert \,\textrm{ord}P\leq i\} \subset \SD_X \qquad (i \in \bbZ)
\end{equation*}
we obtain a filtration $\{F_i\SD_X\}_{i \in \bbZ}$ of $\SD_X$. 
Then by the relation 
$[\partial_{x_i},x_j] = \delta_{ij} \,(1\leq i,j \leq n)$, for the classes 
$\xi_i \coloneqq [\partial_{x_i}] \in F_1\SD_X / F_0\SD_X \,(1\leq i \leq n)$ the sheaf of graded rings 
\begin{equation*}
\gr^F\SD_X \coloneqq \bigoplus_{i\geq0} \Bigl(F_i\SD_X / F_{i-1}\SD_X \Bigr)
\end{equation*}
associated to it is locally isomorphic to $\SO_X[\xi_1,\ldots,\xi_n]$. 
This in particular implies that $\gr^F\SD_X$ is commutative. 
Building upon the coherency of $\SO_X$ proved by K.Oka and using this 
commutative $\gr^F\SD_X$, Kashiwara overcame the difficulty arising from 
the non-commutativity of $\SD_X$ and proved that the 
sheaf $\SD_X$ of (non-commutative) rings is coherent. 
This implies that a (left) $\SD_X$-module $\SM$ is coherent over $\SD_X$ if 
and only if it 
is locally finitely presented over $\SD_X$ i.e. 
for any point $p \in X$ there exist its neighborhood $U \subset X$, 
non-negative integers $N_0,N_1 \in \bbZ_+$ and an exact sequence 
\begin{equation*}
\SD_X^{N_1} \stackrel{\Phi}\longrightarrow \SD_X^{N_0} \longrightarrow \SM \longrightarrow 0
\end{equation*}
of (left) $\SD_X$-modules on $U$. 
As for such $U \subset X$ there exists an $N_1\times N_0$ matrix 
$P = (P_{ij})_{1\leq i\leq N_1,1\leq j\leq N_0} \in M\bigl(N_1,N_0 ; 
\SD_X(U)\bigr)$ of differential operators 
$P_{ij} \in \SD_X(U)$ on $U$ such that the $\SD_X$-linear morphism $\SD_X^{N_1} 
\stackrel{\Phi}\longrightarrow \SD_X^{N_0}$ 
is given by 
\begin{equation*}
\Phi (Q_1,\ldots,Q_{N_1}) = (Q_1,\ldots,Q_{N_1})P \qquad \bigl((Q_1,\ldots,Q_{N_1}) 
 \in \SD_X^{N_1}(U) \bigr), 
\end{equation*}
we obtain an isomorphism $\SM \simeq \SD_X^{N_0} / \SD_X^{N_1}P$ on $U \subset X$. 
In the early 1960's, 
in his faculty talk at the University of Tokyo, 
Mikio Sato proposed a new program to study general systems of 
linear partial differential equations by regarding 
them as coherent $\SD_X$-modules. 
Let us explain his brilliant idea briefly. 
First, note that the sheaf $\SO_X$ is naturally a left $\SD_X$-module. 
Next, apply the left exact functor 
$\shom_{\SD_X}(\,\cdot\,,\SO_X)$ to the exact sequence 
\begin{equation*}
\SD_X^{N_1} \stackrel{\times P}\longrightarrow \SD_X^{N_0} \longrightarrow \SM \longrightarrow 0
\end{equation*}
on $U \ \Bigl(P \in M\bigl( N_1,N_0;\SD_X(U) \bigr) \Bigr)$. 
Then one obtains an exact sequence 
\begin{equation*}
0 \longrightarrow \shom_{\SD_X}(\SM,\SO_X) \longrightarrow \SO_X^{N_0} 
\stackrel{P\times}\longrightarrow \SO_X^{N_1}
\end{equation*}
and hence an isomorphism 
\begin{equation*}
\shom_{\SD_X}(\SM,\SO_X) \simeq \bigl\{ \vec{u} \in \SO_X^{N_0} 
\,\bigm\vert \, P\vec{u} = \vec{0} \bigr\} \qquad \subset \SO_X^{N_0}.
\end{equation*}
We thus find that the sheaf of the holomorphic solutions to the 
differential equation $P\vec{u} = \vec{0}$ is isomorphic to the hom 
sheaf $\shom_{\SD_X}(\SM,\SO_X)$. If $X$ is a complexification of 
a real analytic manifold $M \subset X$, the same is true also after replacing 
$\SO_X$ by any other left $\SD_X|_M$-module of 
classical function spaces such as (complex-valued) $C^\infty$-functions, 
distributions and Sato's hyperfunctions etc. on $M$. 
For a coherent $\SD_X$-module $\SM$ on $X$, then it is natural 
to consider not only the solution sheaf 
$\shom_{\SD_X}(\SM,\SO_X)$ but also the solution complex $\Sol_X(\SM)$ defined as follows. Let
\begin{equation*} 
\cdots \cdots \longrightarrow \SD_X^{N_k} \stackrel{\times P_k}\longrightarrow \SD_X^{N_{k-1}} 
\stackrel{\times P_{k-1}}\longrightarrow \cdots \cdots \longrightarrow 
\SD_X^{N_1} \stackrel{\times P_1}\longrightarrow \SD_X^{N_0} \longrightarrow \SM
\end{equation*}
be a (locally defined) free resolution of $\SM$ over $\SD_X$. Then we call the complex 
\begin{equation*}
0 \longrightarrow \SO_X^{N_0} \stackrel{P_1\times}\longrightarrow \SO_X^{N_1} 
\stackrel{P_2\times}\longrightarrow \cdots \cdots \longrightarrow 
\SO_X^{N_{k-1}} \stackrel{P_k\times}\longrightarrow \SO_X^{N_k} \longrightarrow 
\cdots \cdots 
\end{equation*}
obtained by applying the functor $\shom_{\SD_X}(\,\cdot \,,\SO_X)$ to 
it the solution complex of $\SM$ and denote it by $\Sol_X(\SM)$. 
If we consider it as an object of the derived category $\Db(X)$ 
(for the definition, see e.g. \cite{H-T-T}, \cite{KS2}, 
\cite{KS3}) consisting of 
bounded complexes of $\bbC_X$-modules on $X$, this definition does 
not depend on the choice of the resolution of $\SM$. 
Namely, we set $\Sol_X(\SM) = \rhom_{\SD_X}(\SM,\SO_X)$. 
Obviously it satisfies the condition
\begin{equation*}
H^0\Sol_X(\SM) \simeq \shom_{\SD_X}(\SM,\SO_X).
\end{equation*}
In order to study deeper structures of coherent $\SD_X$-modules, Kashiwara 
defined their characteristic varieties as follows. Let 
$\pi_X \colon T^*X \rightarrow X$ be the (holomorphic) cotangent 
bundle of X. Let $x=(x_1,\ldots,x_n)$ be a local coordinate of $X$ 
and $(x,\xi)=(x_1,\ldots,x_n,\xi_1,\ldots,\xi_n)$ that of $T^*X$ associated to it. 
Then for $P= \sum_{\alpha \in \bbZ_+^n}a_\alpha(x) \partial_x^\alpha \in \SD_X$ of order $m\in \bbZ_+$ 
we define a holomorphic function $\sigma(P)$ of $(x,\xi)$ by 
\begin{equation*}
\sigma(P)(x,\xi) \coloneqq \sum_{\vert \alpha \vert =m} a_\alpha(x)\xi^{\alpha}
\end{equation*}
and call it the principal symbol of $P$. We can easily check that it defines a function on $T^*X$. 
More precisely $\sigma(P)$ is a section of the subsheaf 
\begin{equation*}
\gr^F \SD_X \simeq \SO_X[\xi_1,\ldots,\xi_n] \quad \subset {\pi_X}_*\SO_{T^*X}
\end{equation*}
of the direct image sheaf ${\pi_X}_*\SO_{T^*X}$ of $\SO_{T^*X}$ by $\pi_X$. 
Recall that in classical analysis for a single equation 
$Pu=0 \, (P\in \SD_X)$ (which corresponds to the coherent 
$\SD_X$-module $\SD_X / \SD_X P$) the complex hypersurface 
$\{\sigma(P)=0\} \subset T^*X$ in $T^*X$ is called the characteristic 
variety of the differential operator $P \in \SD_X$ and plays an 
important role in studying the properties of its solutions. However, 
it was not until Kashiwara's master thesis \cite{Kas12} 
that we knew how to define characteristic varieties for general systems 
$P\vec{u}=\vec{0}$ of linear partial differential equations. 
Namely, it was Kashiwara who first set up a solid basis for the 
systematic study of linear partial differential equations. 
For the definition of the characteristic variety of a coherent 
$\SD_X$-module $\SM$ we proceed as follows. 
First, by the (locally defined) surjective morphism $\Psi \colon 
\SD_X^{N_0} \longtwoheadrightarrow\SM$ 
we define a filtration $\{F_i\SM\}_{i\in \bbZ}$ of $\SM$ by
\begin{equation*}
F_i\SM \coloneqq \Psi \bigl((F_i\SD_X)^{N_0} \bigr) \subset \SM \qquad (i\in \bbZ).
\end{equation*}
Then the graded module 
\begin{equation*}
\gr^F \SM \coloneqq \bigoplus_{i\in\bbZ} (F_i\SM / F_{i-1}\SM)
\end{equation*}
associated to it is naturally a $\gr^F\SD_X$-module. Moreover, 
we can show that it is coherent over $\gr^F\SD_X$. 
Then by the ring homomorphisms
\begin{equation*}
\pi_X^{-1}\gr^F\SD_X \longrightarrow \pi_X^{-1}{\pi_X}_*\SO_{T^*X}\longrightarrow \SO_{T^*X}
\end{equation*}
we define a coherent $\SO_{T^*X}$-module $\widetilde{\SM}$ on $T^*X$ by
\begin{equation*}
\widetilde{\SM} \coloneqq \SO_{T^*X} \otimes_{\pi_X^{-1}\gr^F\SD_X}\pi_X^{-1}\gr^F\SM
\end{equation*}
and set $\ch\SM \coloneqq \supp \widetilde{\SM} \subset T^*X$. 
We can show that this definition of the complex analytic 
subset $\ch\SM \subset T^*X$ does not depend on the choice of the 
surjective morphism $\Psi \colon \SD_X^{N_0} 
\longtwoheadrightarrow \SM$ with the help of the theory of 
good filtrations (see \cite{Kas12}). 
We call $\ch\SM \subset T^*X$ thus obtained the characteristic variety of $\SM$. 
As in the case of a single equation 
$Pu=0$, the characteristic variety $\ch\SM \subset T^*X$ is conic i.e. 
invariant under the standard action of 
the multiplicative group $\bbC^*=\bbC \setminus \{0\}$ on $T^*X$. Recall that $T^*X$ is endowed with 
the canonical symplectic structure defined by the non-degenerate holomorphic 2-form 
$\sigma_{T^*X} = \sum_{i=1}^n dx_i\wedge d\xi_i \in \Omega_{T^*X}^2$. 
Then the following theorem obtained by Sato, Kawai and Kashiwara in 
the monumental paper \cite{SKK} 
reveals a striking feature of 
characteristic varieties. 
\begin{theorem}\label{thm:involutive}
The characteristic variety $\ch\SM \subset T^*X$ of a coherent $\SD_X$-module $\SM$ is involutive with respect to 
the canonical symplectic structure of $T^*X$. In particular, for 
any irreducible component $\Lambda$ of $\ch\SM$ we have 
\begin{equation*}
\dim\,\Lambda \geq \frac{1}{2}\dim T^*X =\dim X.
\end{equation*}
\end{theorem}
In the algebraic case, the inequality 
$\dim\, \ch\SM \geq \dim X$ was proved also by 
Bernstein \cite{Ber1} and we call it Bernstein's inequality. 
Later, a purely algebraic proof of Theorem \ref{thm:involutive} was 
given by Gabber \cite{Ga}. 
For a geometric proof of this ``involutivity theorem" 
see Kashiwara-Schapira \cite[Theorem 6.5.4]{KS2}, 
in which to our surprise the involutivity 
of the micro-supports of general sheaves is proved. Among coherent 
$\SD_X$-modules, the following ones are of particular importance.

\begin{definition}
We say that a coherent $\SD_X$-module $\SM$ is holonomic if its 
characteristic variety $\ch\SM \subset T^*X$ is Lagrangian 
i.e. $\dim \, \ch\SM \leq \dim X$.
\end{definition}

Let $\Mod_\mathrm{coh}(\SD_X)$ be the abelian category of coherent $\SD_X$-modules on $X$ and 
$\Mod_\mathrm{h}(\SD_X) \subset \Mod_\mathrm{coh}(\SD_X)$ its full 
subcategory consisting of holonomic $\SD_X$-modules. 
Then it turns out that the category $\Mod_\mathrm{h}(\SD_X)$ is again 
abelian. It was Mikio Sato who first noticed 
the importance of holonomic $\SD_X$-modules and suggested Kashiwara 
and Kawai to study them very precisely. 
Indeed, as in the (quite fruitful) case of ordinary differential 
equations on the complex plane $\bbC$ in classical analysis, the holomorphic solution complex 
$\Sol_X(\SM)= \rhom_{\SD_X}(\SM,\SO_X) \in \Db(X)$ is ``finite dimensional" 
in the following sense. Recall that for 
a complex manifold $Z$ a $\bbC_Z$-module $L$ is a local system if it is locally 
free of finite rank over $\bbC_Z$. For the complex manifold $X$, then we 
say that a $\bbC_X$-module $K$ is (complex) constructible if there exists a (complex analytic) 
stratification $X= \bigsqcup_{\alpha\in A}X_\alpha$ of $X$ such that 
$K\vert_{X_\alpha}$ is a local system for any $\alpha\in A$. 
More generally, we say that an object $F\in \Db(X)$ of the 
derived category of $\Db(X)$ of $\bbC_X$-modules is constructible if 
its cohomology sheaves are constructible in the above sense. The following theorem, which is now called 
Kashiwara's constructibility theorem, was proved in his Ph.D. thesis \cite{Kas3}.

\begin{theorem}
For a holonomic $\SD_X$-module $\SM$ its solution 
complex $\Sol_X(\SM) = \rhom_{\SD_X}\\(\SM,\SO_X) \in \Db(X)$ is constructible.
\end{theorem}

For a holonomic $\SD_X$-module $\SM$ we set 
\begin{equation*}
\DR_X(\SM) \coloneqq \rhom_{\SD_X}(\SO_X,\SM)[\dim X]
\end{equation*}
and call it the de Rham complex of $\SM$. In \cite{Kas3} Kashiwara also 
proved that $\DR_X(\SM) \in \Db(X)$ is constructible. 
Let $\Dbc(X) \subset \Db(X)$ be the full subcategory of $\Db(X)$ 
consisting of constructible objects. 
Recall that the Verdier dual functor $\mathrm{D}_X \colon \Db(X) \longrightarrow \Db(X)$ defined by 
\begin{equation*}
\mathrm{D}_X(F) \coloneqq \rhom_{\bbC_X}(F,\bbC_X[2\dim X]) \qquad \bigl( F \in \Db(X) \bigr)
\end{equation*}
preserves the constructibility. Then the relation between $\Sol_X(\SM)$ 
and $\DR_X(\SM)$ is explained by the isomorphism 
\begin{equation*}
\DR_X(\SM) \simeq \mathrm{D}_X \bigl( \Sol_X(\SM)[\dim X] \bigr)
\end{equation*}
proved also in Kashiwara \cite{Kas3}. Although the notion of perverse 
sheaves appeared only after the paper \cite{BBD} of Beilinson, 
Bernstein and Deligne in the 1980's, in the Ph.D. thesis 
\cite{Kas3} Kashiwara already proved that 
for a holonomic $\SD_X$-module $\SM$ the constructible sheaves 
$\Sol_X(\SM)[\dim X]$, $\DR_X(\SM) \in \Dbc(X)$ 
satisfy the conditions of perversity defined as follows.

\begin{definition}
We say that $F \in \Dbc(X)$ is a perverse sheaf if it satisfies the following two conditions.
\begin{enumerate}
\item[(i)] $\dim \,(\supp H^j F) \leq -j \qquad (j\in \bbZ)$,
\item[(ii)] $\dim \, \bigl( \supp H^j \mathrm{D}_X(F) \bigr) \leq -j \qquad (j\in \bbZ)$.
\end{enumerate}
\end{definition}
Denote by $\Perv(\bbC_X) \subset \Dbc(X)$ the full subcategory 
of $\Dbc(X)$ consisting of perverse sheaves. Then it turns out 
that the category $\Perv(\bbC_X)$ is abelian. Consequently, we obtain a functor 
\begin{equation*}
\DR_X(\, \cdot \,) \colon \Mod_\mathrm{h}(\SD_X) \longrightarrow\Perv(\bbC_X).
\end{equation*}
However it is not an equivalence of categories in general except for 
the trivial case of $\dim X=0$. A remedy for 
this problem is to construct a suitable subcategory of $\Mod_\mathrm{h}(\SD_X)$ to have an equivalence. 
Generalizing the classical notion of regular singular points of 
ODEs, in \cite{KO} Kashiwara and Oshima defined the regularity of holonomic $\SD$-modules. 
Then in \cite{KK1} Kashiwara and Kawai developed the theory of 
regular holonomic $\SD$-modules and proved their basic properties. 
Finally in \cite{Kas6} and \cite{Kas9} 
Kashiwara established the following 
Riemann-Hilbert correspondence. Let $\Mod_\mathrm{rh}(\SD_X)$ be the 
full subcategory of $\Mod_\mathrm{h}(\SD_X)$ consisting of regular holonomic $\SD_X$-modules. 
We denote by $\Mod_\mathrm{rh}(\SD_X)^{\mathrm{op}}$ its opposite category.
\begin{theorem}
The functors 
\begin{align*}
\DR_X(\, \cdot \,) \colon & \Mod_\mathrm{rh}(\SD_X) \longrightarrow\Perv(\bbC_X) \\
\Sol_X(\, \cdot \,)[\dim X] \colon & \Mod_\mathrm{rh}(\SD_X)^{\mathrm{op}}\longrightarrow\Perv(\bbC_X)
\end{align*}
induce equivalences of categories.
\end{theorem}

Let $\Db(\SD_X)$ be the derived category of bounded complexes of 
(left) $\SD_X$-modules and $\Dbrh(\SD_X)$ its full subcategory consisting 
of bounded complexes of $\SD_X$-modules whose cohomology sheaves are regular holonomic.
Then, more generally in \cite{Kas9} Kashiwara proved also that there 
exist equivalences of categories 
\begin{equation*}
\DR_X(\,\cdot\,)\colon \Dbrh(\SD_X)\simto\Dbc(X), \quad
\Sol_X(\,\cdot\,)\colon \Dbrh(\SD_X)^\op\simto\Dbc(X).
\end{equation*}
This gives an ultimate solution to (the higher-dimensional generalization of) 
the 21st problem of Hilbert, raised in the early 20th century.
Recall that for regular meromorphic connections $\SM\in\Modrh(\SD_X)$ 
on $X$ a satisfactory answer to it was given previously by Deligne \cite{De1}. 
In order to prove the Riemann-Hilbert correspondence, Kashiwara 
constructed an inverse $\RH_X(\,\cdot\,)\colon\Dbc(X)\longrightarrow\Dbrh(\SD_X)^\op$ of
the solution functor $\Sol_X(\,\cdot\,)\colon\Dbrh(\SD_X)^\op 
\longrightarrow\Dbc(X)$ explicitly. 
His proof was completed around 1980 (see \cite{Kas6}). 
Then the full proof was published in \cite{Kas9}. 
Later, a different proof was given by Mebkhout in 
\cite{Me2}. The Riemann-Hilbert correspondence has also a nice 
functoriality with respect to morphisms of complex 
manifolds as follows. First, note that for a morphism 
$f\colon Y \longrightarrow X$ of complex manifolds 
the right $f^{-1}\SD_X$-module 
\begin{equation*}
\SD_{Y \rightarrow X} 
\coloneqq \SO_Y \otimes_{f^{-1}\SO_X} f^{-1}\SD_X
\end{equation*}
is naturally endowed with a structure of a left 
$\SD_Y$-module, commuting with the action of $f^{-1}\SD_X$. 
Namely $\SD_{Y \rightarrow X}$ is a 
$(\SD_Y, f^{-1}\SD_X)$-bimodule and for 
$\SM \in \Db(\SD_X)$ we define its inverse image 
$\mathrm{D} f^\ast \SM \in \Db(\SD_Y)$ by 
\begin{equation*}
\mathrm{D} f^\ast \SM
\coloneqq \SD_{Y \rightarrow X} 
\overset{L}\otimes_{f^{-1}\SD_X} f^{-1} \SM 
\qquad \in \Db(\SD_Y). 
\end{equation*}
Since this operation preserves the regular holonomicity, 
we obtain a functor 
\begin{equation*}
\mathrm{D} f^\ast ( \,\cdot\, ): 
\Dbrh(\SD_X) \longrightarrow \Dbrh(\SD_Y)
\end{equation*}
(see \cite{Kas5} and \cite{KK1}). 
Moreover, for any $\SM \in \Dbrh(\SD_X)$ there exists 
an isomorphism 
\begin{equation*}
f^{-1} \Sol_X(\SM) \simeq \Sol_Y(\mathrm{D} f^\ast \SM). 
\end{equation*}
Similarly, if $f$ is proper for $\SN \in \Dbrh(\SD_Y)$ 
we can define its direct image 
$\mathrm{D} f_\ast \SN \in \Dbrh(\SD_X)$ and obtain 
an isomorphism 
\begin{equation*}
\mathrm{R} f_\ast \Sol_Y(\SN ) [\dim Y]
\simeq 
\Sol_X(\mathrm{D} f_\ast \SN ) [\dim X]. 
\end{equation*}
Via the Riemann-Hilbert correspondence, problems in analysis ($\SD$-modules) 
can be interpreted to those in geometry (constructible sheaves). 
In this way, many problems were settled and beautiful theories 
have emerged, advancing the state of science of mathematics tremendously.
\footnote{For instance, motivated by the Riemann-Hilbert correspondence 
and the theory of perverse sheaves, in \cite{GM1} Goresky and 
MacPherson introduced a new notion of intersection 
cohomology groups. See \cite{H-T-T}, \cite{Kir} and 
\cite{Mx} for the details.} 
Among other things, in Brylinski-Kashiwara \cite{BK} 
and Beilinson-Bernstein \cite{BB} the Kazhdan-Lusztig conjecture 
in \cite{KL1} on the character formula for some standard representations of 
semisimple Lie algebras was proved. 
For this purpose, the authors of \cite{BK} and \cite{BB} constructed $\SD$-modules 
on flag manifolds and used the perverse sheaves corresponding to 
them by the Riemann-Hilbert correspondence to relate the characters of 
the representations to the Kazhdan-Lusztig polynomials. 
Later, in \cite{KT1}, \cite{KT2} and \cite{KT3} 
Kashiwara and Tanisaki settled also similar conjectures 
on affine Lie algebras. For this purpose, they 
used $\SD$-modules on infinite dimensional flag schemes. 
For the further development of the theory of the Riemann-Hilbert correspondence 
for regular holonomic $\SD$-modules, 
see e.g. \cite{An1}, \cite{An2}, \cite{Ber3}, \cite{ALD}, 
\cite{Col}, \cite{F-MF-S}, \cite{H-T-T}, 
\cite{KS96}, \cite{KS16}, \cite{Neto}, \cite{Was}. 
For the computational aspect of the theory of $\SD$-modules, 
see e.g. \cite{SST}.  

With Kashiwara's constructibility theorem at hand, it is also natural 
to ask if for a holonomic $\SD_X$-module $\SM$ its local Euler-Poincar{\'e} indices
\begin{equation*}
\chi_x(\Sol_X(\SM))\coloneq\sum_{j\in\bbZ}(-1)^j
\dim_\bbC H^j\Sol_X(\SM)_x \quad (x\in X)
\end{equation*}
can be expressed in terms of $\SM$.
Surprisingly, a beautiful answer to this problem was given also in Kashiwara's Ph.D. 
thesis as follows (for the details, see \cite{Kas2} and \cite{Kas7}).  
Let $\ch\SM=\bigcup_{i\in I}\Lambda_i$ be the irreducible 
decomposition of the characteristic variety $\ch\SM=\supp\tl{\SM}\subset T^\ast X$ 
of $\SM$ and for $i\in I$ denote by $m_i>0$ the multiplicity of the coherent 
$\SO_{T^\ast X}$-modules $\tl{\SM}$ along the irreducible component $\Lambda_i\subset\supp\tl{\SM}$.
Then we obtain a Lagrangian cycle
\begin{equation*}
\CCyc(\SM)\coloneq\sum_{i\in I}m_i \cdot [\Lambda_i]
\end{equation*}
in $T^\ast X$ and call it the characteristic cycle of $\SM$.
For $i\in I$ let us consider the (irreducible) analytic subset 
$Y_i\coloneq\pi_X(\Lambda_i)\subset X$ of $X$ and its regular (smooth) 
part $Z_i\coloneq (Y_i)_\reg\subset Y_i$.
Then one can show that the closure $\overline{T_{Z_i}^\ast X}$ of the 
conormal bundle $T_{Z_i}^\ast X\subset T^\ast X$ of $Z_i\subset X$ in $T^\ast X$ coincides with $\Lambda_i$.
In this situation, using stratifications of $Y_i\subset X$ Kashiwara 
defined a $\bbZ$-valued constructible function 
$\Eu_{Y_i}\colon Y_i\longrightarrow \bbZ$ on $Y_i$ such that 
$\Eu_{Y_i}(x)=1$ for any $x\in Z_i=(Y_i)_\reg\subset Y_i$
and proved the following theorem, which is now called Kashiwara's index theorem.
Here we say that a function $\varphi$ on an analytic set $W$
is constructible if there exists a stratification $\calS$ of $W$ 
such that $\varphi$ is 
constant on each stratum $S\subset W$ in $\calS$.

\begin{theorem}
In the situation as above, for any $x\in X$ we have 
\begin{equation*}
\chi_x(\Sol_X(\SM))=\sum_{i\colon x\in Y_i}(-1)^{\codim_X Y_i}
\cdot m_i \cdot \Eu_{Y_i}(x).
\end{equation*}
\end{theorem}

We call the $\bbZ$-valued constructible function $\Eu_{Y_i}\colon Y_i\longrightarrow \bbZ$ 
on $Y_i\subset X$ the Euler obstruction of $Y_i$.
This function plays an important role in singularity theory and 
was found also by MacPherson in \cite{Mac} independently.
He used it to prove the existence of characteristic classes of singular 
varieties that are functorial with respect to their morphisms and 
thereby proved a conjecture of Grothendieck and Deligne.
We now call them the Chern-Schwartz-MacPherson classes.
Later, by Sabbah \cite{Sab1} and Kennedy \cite{Ken} their 
functoriality was more neatly explained by that of the 
characteristic cycles of holonomic $\SD$-modules.
For the further development of Euler obstructions and 
Chern-Schwartz-MacPherson classes, see e.g. \cite{B-S-Y}, \cite{Er},  
\cite{G-G-R}, \cite{Go}, \cite{M-T-3}, \cite{P-P}, \cite{Schu1}, 
\cite{Schu2} etc.

\section{Microlocal analysis}
The idea of symplectic geometry and microlocal analysis goes back to 
the age of Hamilton and Jacobi. 
It is Hamilton who first considered ordinary differential equations 
on the cotangent bundle $T^\ast\mathbb{R}^n$ of $\mathbb{R}^n$ in order to
study those on $\mathbb{R}^n$ more elegantly. 
Such systems of ODEs are called Hamilton's canonical equations 
and are basic in analytical mechanics and classical analysis. 
Jacobi used Hamilton's ODEs to construct solutions of partial
differential equations of first order. 
We call it the method of Hamilton and Jacobi. 
The idea of microlocal analysis of Mikio Sato is to see 
the singularities of functions on $\mathbb{R}^n$ locally 
at each point of the cotangent bundle $T^\ast\mathbb{R}^n$ and 
extend the above mentioned beautiful results in classical analysis 
to any system of linear partial differential equations, 
at least on the singularities of their solutions. 
In collaboration with Kawai and Kashiwara, his dream was 
realized in the so-called SKK paper \cite{SKK}, 
revealing numerous astonishing results. 
In order to explain some of their great achievements in it, 
let $M$ be a real analytic manifold and $X$ its complexification. 
We denote by $T_M^*X\simeq \sqrt{-1} T^*M$ the conormal bundle of $M$ in $X$. 
By removing the zero section $T_M^\ast X\cap T_X^\ast X\simeq M\subset T_M^\ast X$ 
from it we set $\mathring{T}_M^\ast X\coloneq T_M^\ast X\setminus T_X^\ast X$. 
Let $\Sa_M \coloneq \SO_X|_M$ (resp. $\mathcal{B}_M$) be 
the sheaf of real analytic functions (resp. Sato's hyperfunctions) on $M$. 
Then in 1969, Sato introduced the sheaf $\SC_M$ of microfunctions on 
$\mathring{T}_M^\ast X$ \footnote{In fact, in \cite{SKK} 
the sheaf of microfunctions was defined on the 
sphere bundle $S_M^\ast X= \mathring{T}_M^\ast X 
/ \bbR_{>0}$.} which fits into the exact sequence 
\begin{equation*}
0\longrightarrow \Sa_M \longrightarrow \mathcal{B}_M \longrightarrow 
\mathring{\pi}_\ast \SC_M \longrightarrow 0
\end{equation*}
where $\mathring{\pi} \colon \mathring{T}_M^\ast X \to M$ is the projection. 
For this purpose, he defined the microlocalization functor 
\begin{equation*}
\mu_M ( \cdot ): \Db(X) \longrightarrow \Db(T^*_MX)
\end{equation*}
along $M \subset X$, which is now a basic tool in sheaf theory. 
Moreover, Kashiwara proved that the sheaf $\SC_M$ is conically flabby. 
For a hyperfunction $u\in \mathcal{B}_M(U)$ on an open subset $U\subset M$, 
by the morphism $\mathcal{B}_M \longrightarrow \mathring{\pi}_\ast \SC_M$ 
we thus obtain a section $v\in \SC_M(\mathring{\pi}^{-1}(U))$ of $\SC_M$ on 
$\mathring{\pi}^{-1}(U) \subset \mathring{T}_M^\ast X$. 
Its support is nothing but the analytic wave front set 
${\rm WF}(u)$ of $u\in \mathcal{B}_M(U)$. 
Then by the exact sequence 
\begin{equation*}
0\longrightarrow \Sa_M (U) \longrightarrow \mathcal{B}_M (U) \longrightarrow 
\SC_M(\mathring{\pi}^{-1}(U))
\end{equation*}
we see that $u\in \Sa_M(U) \iff {\rm WF}(u)= \emptyset$. This 
shows that the microfunction $v\in \SC_M(\mathring{\pi}^{-1}(U))$ can be 
regarded as the ``singularity" of $u\in \mathcal{B}_M (U)$. 
In \cite{SKK} the authors also introduced the sheaf $\SE_X$ of 
(non-commutative) rings of microdifferential operators on 
the complex cotangent bundle 
$T^\ast X$ of $X$ such that for a coherent $\SD_X$-module $\SM$ we have 
\begin{equation*}
\ch \SM = \supp \left( \SE_X \otimes_{\pi_X^{-1} \SD_X} \pi_X^{-1} \SM \right).
\end{equation*}
Moreover they showed that its restriction $\SE_X|_{\mathring{T}_M^\ast X}$ 
to $\mathring{T}_M^\ast X\subset T^\ast X$ acts on the sheaf 
$\SC_M$ of microfunctions. 
For a coherent $\SD_X$-module $\SM$ 
we define a coherent $\SE_X$-module $\SE \SM$ on $T^\ast X$ by
\begin{equation*}
\SE\SM\coloneq \SE_X\otimes_{\pi_X^{-1}\SD_X}\pi_X^{-1}\SM.
\end{equation*}
Then we obtain an isomorphism 
\begin{equation*}
\rhom_{\SE_X}(\SE\SM,\SC_M)\simeq\rhom_{\pi_X^{-1}\SD_X}(\pi_X^{-1}\SM,\SC_M),
\end{equation*}
where the restrictions $(\,\cdot\,)|_{\mathring{T}_M^\ast X}$ to 
$\mathring{T}_M^\ast X\subset T^\ast X$ are omitted. 
Hence, if the coherent $\SD_X$-module $\SM$ is elliptic i.e. 
if it satisfies the condition
\begin{equation*}
 \supp ( \SE \SM ) \cap \mathring{T}_M^\ast X = 
\ch \SM \cap \mathring{T}_M^\ast X=\varnothing
\end{equation*}
we can easily see that there exists an isomorphism
\begin{equation*}
\rhom_{\SD_{X}|_M} (\SM|_M, \Sa_M) \overset{\sim}{\longrightarrow} 
\rhom_{\SD_{X}|_M} (\SM|_M, \mathcal{B}_M).
\end{equation*}
Namely, in this case, all the hyperfunction solutions to $\SM$ are real analytic. 
Without the sheaf $\SC_M$, it would be impossible to prove 
such a very general regularity theorem for systems of linear PDEs. 
We call it Sato's fundamental theorem. 
More generally, if a hyperfunction $u\in \mathcal{B}_M(U)$ on an open 
subset $U\subset M$ is a solution of a coherent $\SD_X$-module $\SM$, its 
analytic wave front set ${\rm WF}(u)$ is contained 
in the set $\ch \SM \cap \mathring{T}_M^\ast X$. 
Then the next task of the authors of \cite{SKK} was 
to study the structures of the microfunction solutions 
\begin{equation*}
\rhom_{\SE_X}(\SE\SM,\SC_M)\simeq\rhom_{\pi_X^{-1}\SD_X}(\pi_X^{-1}\SM,\SC_M)
\end{equation*}
of general coherent $\SD_X$-modules $\SM$. 
First in \cite{SKK}, they developed a theory of quantized contact transformations
of coherent $\SE_X$-modules i.e. 
microdifferential systems and proved that by a (locally defined) 
homogeneous symplectic transform of $T^\ast X$ 
(and microdifferential operators of infinite order) 
the coherent $\SE_X$-module 
$\SE \SM$ can be locally transformed to a partial de Rham system $\partial_{x_1} 
u=\cdots =\partial_{x_d} u=0$
in general. 
Furthermore, they showed that by a (locally defined) 
``real" homogeneous symplectic transform 
of $T^\ast X$ preserving $T_M^\ast X\subset T^\ast X$, 
$\SE \SM$ is locally 
transformed to a mixture of partial de Rham, partial Cauchy-Riemann 
and Lewy-Mizohata systems in general. 
It is surprising that we have such a general structure theorem 
for generic microdifferential systems. 
From this, we see in particular that 
if $V\coloneq \ch \SM \cap \mathring{T}_M^\ast X\subset \mathring{T}_M^\ast X$ 
is an involutive submanifold of $T_M^\ast X\simeq \sqrt{-1} T^\ast M$ 
and $\SM$ has real simple characteristic along $V$ at a point of $V$ 
then the microfunction solutions 
of $\SM$ are constant along the bicharacteristic leaves of $V$, 
which means that the singularities of the solutions of $\SM$ propagate along them.
In the special case where $\SM$ is the classical wave equation, 
this gives the first mathematical 
explanation to the particle nature of waves, 
which is a traditional guiding principle in physics. 
Note that inspired by Sato's ideas, in \cite{H} and \cite{D-H} 
H\"ormander introduced the wave front sets of distributions 
and developed a theory of Fourier integral operators, which is 
parallel to the one of QCT of SKK, in collaboration with Duistermaat.

\section{Microlocal study of sheaves}
After the microlocal study of hyperbolic $\SD$-modules in \cite{KS79},
Kashiwara and Schapira gradually noticed that in order to obtain
the main results on the microfunction and hyperfunction solutions 
to a coherent $\SD$-module the only thing that they have to keep in mind
is  the set of the codirections (in the cotangent bundle) in which
the cohomology groups of its solution complex do not propagate.
First, by abstracting the notion of characteristic varieties 
of coherent $\SD$-modules,
for any bounded complex of sheaves $F$ on a manifold they succeeded in defining
a closed conic subset in its cotangent bundle and called it the micro-support
of $F$. Note that this definition was motivated also by 
the previous works by Bony-Schapira \cite{BS} 
and Zerner \cite{Z}.  
Then in \cite{KS1} and \cite{KS2}, they studied the functorial properties of the 
micro-supports of sheaves on manifolds and developed the whole  theory in a short period.
It turned out that this theory is extremely useful in many branches of contemporary
mathematics.
In this section, without getting into the details of the theory, we shall 
introduce the main results in \cite{KS1} and \cite{KS2} and explain related results.
For the sake of simplicity, 
in this section we assume that $\mathbf{k}$ is a commutative field
of characteristic zero and for a (real and $C^{\infty}$) manifold $X$
we denote by $\Db(X)$ the derived category of bounded complexes of
sheaves of $\mathbf{k}_X$-modules on it.
Let $\pi_X \colon T^\ast X \longrightarrow X$ be the projection.

\begin{definition}
For an object $F\in \Db(X)$ we define a subset $\MS (F)\subset T^\ast X$
by:\par
\noindent For a point $p\in T^\ast X$ we have
\begin{align*}
p \notin \MS(F)
\iff
\left\{
\begin{array}{l}
\text{There exists an open neighborhood $U$ of $p$ in $T^\ast X$ such that} \\
\text{for any point $x\in \pi_X(U) \subset X$ and any $C^{\infty}$-function
$\phi$ defined} \\
\text{on a neighborhood of $x$
in $X$ satisfying the conditions}\\
\begin{aligned}
\phi (x)=0, \quad \mathrm{d} \phi (x)\in U 
\end{aligned}\\
\text{we have the vanishing}\\
\begin{aligned}
\mathrm{R} \Gamma_{\{\phi \geq 0\}} (F)_x \simto 0.
\end{aligned}
\end{array}
\right.
\end{align*}
We call $\MS (F)\subset T^\ast X$ the micro-support of $F$.
\end{definition}

By this definition, it is clear that $\MS(F) \subset T^\ast X$ is a closed subset
of $T^\ast X$ and conic i.e. invariant by the standard action of the multiplicative
group $\mathbb{R}_{>0} \coloneq \{ t\in \mathbb{R} \mid t>0\}$ on $T^\ast X$.
For a closed  subset $Z\subset X$ of a manifold $X$ and the inclusion map
$i_Z \colon Z \longhookrightarrow X$ we denote 
$\mathrm{R} i_Z{}_\ast (\mathbf{k}_Z)\simeq i_Z{}_\ast (\mathbf{k}_Z) \in \Db(X)$ 
simply by $\mathbf{k}_Z \in \Db(X)$.
\begin{example}
For the constant sheaf $\mathbf{k}_X \in \Db(X)$ on $X$, 
its micro-support $\MS(\mathbf{k}_X) \subset T^\ast X$ is equal to 
the zero section $T_X^\ast X \simeq X \subset T^\ast X$ of $T^\ast X$.
More generally, for a closed submanifold $M\subset X$ of $X$ and the sheaf 
$\mathbf{k}_M\in \Db(X)$ on $X$ associated to it, its micro-support 
$\MS(\mathbf{k}_M) \subset T^\ast X$ is equal to the conormal bundle 
$T_M^\ast X\subset T^\ast X$ of $M$ in $X$.
\end{example}
With micro-supports of sheaves at hand, we can extend various results in
Morse theory. 
For instance, together with Kashiwara's 
non-characteristic deformation lemma (see \cite[Proposition 2.7.2]{KS2}), 
for $F\in \Db(X)$ and two open subsets $U, V\subset X$ in $X$ such that $V\subset U$,
under some suitable conditions on $\MS(F)\subset T^\ast X$ we obtain an 
isomorphism
\begin{equation*}
\rsect(U;F) \simto \rsect(V,F).
\end{equation*}
Note that in the theory of 
stratified Morse theory developed by Goresky and MacPherson 
\cite{GM88} similar results were obtained for complex constructible sheaves. 
The micro-supports of sheaves enjoy several nice functorial properties 
with respect to morphisms of manifolds.
For instance, if $f\colon Y\longrightarrow X$ is a morphism of manifolds 
and we take an object $G\in \Db(Y)$ such that  $f$ is proper on its support
$\supp G\subset Y$, then the micro-support 
$\MS \left( \mathrm{R}f_\ast G\right)\subset T^\ast X$ of the direct image
$\mathrm{R}f_\ast G\in \Db(X)$ is estimated from above as follows.
Let
\begin{equation*}
T^\ast Y \underset{\rho_f}{\longleftarrow} Y\times_X T^\ast X \underset{\varpi_f}{\longrightarrow} T^\ast X
\end{equation*}
be the natural morphisms induced by $f\colon Y\longrightarrow X$.
Then we have 
\begin{equation*}
\MS(\mathrm{R}f_\ast G)\subset \varpi_f \rho_f^{-1} \left( \MS(G) \right) .
\end{equation*}
The theory of micro-supports plays a powerful role also in studying morphisms 
in derived categories.
Indeed, many natural morphisms in them are not isomorphisms.
But with micro-supports, we obtain conditions under which they are actually isomorphisms.
For instance, for a morphism $f\colon Y \longrightarrow X$ of manifolds and $F\in \Db(X)$,
Kashiwara and Schapira proved in \cite[Corollary 6.4.4]{KS2}  
that if $f$ is non-characteristic for
$\MS(F)\subset T^\ast X$ (see \cite[Definition 6.2.7]{KS2} for the definition) 
the natural morphism
\begin{equation*}
f^{-1} F\otimes_{\mathbf{k}_Y} f^{!}\mathbf{k}_X \longrightarrow f^{!}F
\end{equation*}
in the derived category $\Db(Y)$ is an isomorphism. 
As we see in the theorem below, the notion of micro-supports is 
a generalization of that of characteristic varieties.
If $X$ is a complex manifold, then we denote by $X_{\mathbb{R}}$ 
(resp. $(T^*X)_{\mathbb{R}}$) the real analytic
manifold underlying $X$ (resp. the complex cotangent bundle 
$T^\ast X$ of $X$) and identify the cotangent bundle $T^\ast X_{\mathbb{R}}$ 
of $X_{\mathbb{R}}$ with $(T^*X)_{\mathbb{R}}$ as follows. 
For a covector $\mathrm{d} \phi(x) \in T_x^\ast X_{\mathbb{R}}$ at a point
$x\in X_{\mathbb{R}}$ defined by a real-valued $C^{\infty}$-function $\phi$
on a neighborhood of $x$ in $X_{\mathbb{R}}$, by the decomposition
$\mathrm{d} \phi(x) = \partial\phi(x) +\bar{\partial}\phi(x)$
we obtain an element $\partial\phi(x)\in T_x^\ast X$ of $T_x^\ast X$.
Then we see that there exists an isomorphism $T^\ast X_{\mathbb{R}}
\simto (T^\ast X)_{\mathbb{R}}$ (see \cite[Section 11.1]{KS2} 
for a more intrinsic construction of this isomorphism) 
and have the following result (see \cite[Theorem 10.1.1]{KS1}). 

\begin{theorem}[Kashiwara-Schapira \cite{KS1}] 
Let $X$ be a complex manifold and $\SM$ a coherent $\SD_X$-module on it. 
Then for its solution complex $\rhom_{\SD_X} (\SM, \SO_X) \in \Db(X_{\mathbb{R}})$,
under the identification 
$T^\ast X_{\mathbb{R}} \simto (T^\ast X)_{\mathbb{R}}$ as above
we have an equality 
\begin{equation*}
\MS(\rhom_{\SD_X} (\SM, \SO_X)) =\ch \SM. 
\end{equation*}
\end{theorem}
With this theorem at hand, together with the microlocalization functors 
$\mu_M \colon \Db(X) \longrightarrow \Db(T_M^\ast X)$ for submanifolds 
$M\subset X$,
we can recover and even extend various results in microlocal analysis. 
For an interesting generalization of the notion of micro-supports, 
see \cite{K-MF-S}. 
Until the end of this section, we assume that all the manifolds are real analytic.
Then in \cite{KS1} and \cite{KS2} Kashiwara and Schapira introduced the notion
of $\mathbb{R}$-constructible sheaves on such a manifold $X$ by considering
the subanalytic stratifications of $X$ and developed the general theory for them.
Let $\Dbrc(X)\subset \Db(X)$ be the full subcategory of $\Db(X)$
consisting of objects in $\Db(X)$ whose cohomology sheaves are $\mathbb{R}$-constructible.
Then this category is stable by various functors in the derived categories.
For instance, for a morphism $f\colon Y\longrightarrow X$ of real analytic
manifolds and $F\in \Dbrc(X)$ 
we have $f^{-1}F, f^!F\in \Dbrc(Y)$. 
For $G\in \Dbrc(Y)$, if moreover $f$ is proper on the support
$\supp G \subset Y$ of $G$, 
then we have $\rmR f_\ast G, \rmR f_! G\in \Dbrc(X)$.
From this nice property of $\mathbb{R}$-constructible sheaves,
it is natural to study the integral transforms for them 
defined as follows. 
Let $X$ and $Y$ be manifolds and 
\begin{equation*}
X\underset{q_1}{\longleftarrow} X\times Y \underset{q_2}{\longrightarrow} Y
\end{equation*}
the projections.
Assume that $Y$ is compact.
Then for any $\mathbb{R}$-constructible $G\in \Dbrc(Y)$ and 
$K\in \Dbrc(X\times Y)$ the object
\begin{equation*}
\Phi_K(G) \coloneq \mathrm{R}q_1{}_\ast \left(
 K\otimes_{\mathbf{k}_{X\times Y}}^{\mathrm{L}} q_2^{-1}G\right) \simto 
\mathrm{R}q_1{}_\ast \left( K\otimes_{\mathbf{k}_{X\times Y}} q_2^{-1}G\right) \quad \in \Db(X)
\end{equation*}
on $X$ is also $\mathbb{R}$-constructible i.e. $\Phi_K(G)\in \Dbrc(X)$.
We call it the integral transform of $G$ by the kernel $K$.
Although in general it is very hard to calculate $F\coloneq \Phi_K(G)$ explicitly,
by the methods developed by Kashiwara and Schapira in \cite{KS1} and \cite{KS2} 
under some geometric conditions we can determine its microlocal types defined as follows. 
First, as $\MS(F) \subset T^\ast X$ is a conic 
subanalytic Lagrangian subset of $T^\ast X$, 
there exists a subanalytic Whitney stratification 
$X=\bigsqcup_{\alpha \in A} X_\alpha$ of $X$ 
such that $\MS(F) \subset \bigsqcup_{\alpha \in A} T_{X_\alpha}^\ast X$.
This implies that for a generic (smooth) point $p\in \MS(F)$ of $\MS(F)$ 
there exists a submanifold $M\subset X$ such that $\MS(F) =T_M^\ast X$ 
on a neighborhood of $p$ in $T^\ast X$.
Let $\Db(X ;p)$ be the localization of $\Db(X)$ at $p\in \MS(F)$
(see \cite[Section 6.1]{KS2} for the definition).
Let pt be the one point space.
Then by \cite[Propositon 6.6.1]{KS2} there exists an object $L\in \Dbrc(\mathrm{pt})$
such that for the object $L_M \coloneq L_X \otimes_{\mathbf{k}_X} \mathbf{k}_M 
\in \Dbrc(X)$ we have an isomorphism $F \simeq L_M$
in $\Db(X ;p)$.
Let us call $L\in \Dbrc(\mathrm{pt})$ the microlocal type of $F\in \Dbrc(X)$ at $p\in \MS(F)$. 
\footnote{This naive definition is 
different from the one in \cite[Section 7.5]{KS2} 
by some shift.} 
Using this new idea of Kashiwara and Schapira, many mathematicians studied
integral transforms of constructible sheaves and $\SD$-modules in various
situations (see e.g. \cite{Br}, \cite{DE}, \cite{D-S}, \cite{DS2}, 
\cite{Er}, \cite{Mar}, \cite{Mar-Tani}, \cite{ya} etc.). 

Finally, to end this section, we shall introduce 
Kashiwara's microlocal index theorem for 
$\mathbb{R}$-constructible sheaves, which is a generalization of the one
for the solution complexes to holonomic $\SD$-modules proved by
Dubson in \cite{Du}. 
Let $X$ be a manifold and $\orsh_X$ the orientation sheaf on it.
Then by the projection $\pi_X \colon T^\ast X\longrightarrow X$ we set
\begin{equation*}
\mathcal{L}_X \coloneq \varinjlim_{\Lambda} 
H_\Lambda^{\mathrm{dim}X} (\pi_X^{-1} \orsh_X),
\end{equation*}
where $\Lambda$ ranges through the family of all closed conic subanalytic
isotropic subsets of $T^\ast X$. 
Then a section $s\in \Gamma(X; \mathcal{L}_X)$ of the sheaf $\mathcal{L}_X$ 
is considered as a (conic and subanalytic) Lagrangian cycle in $T^\ast X$ 
(twisted by the rank one local system $\pi_X^{-1} \orsh_X$ over 
$\mathbf{k}_{T^\ast X}$) and in \cite{Kas10} 
for an $\mathbb{R}$-constructible sheaf 
$F\in \Dbrc(X)$ on $X$ Kashiwara defined such a section 
$\CCyc(F) \in \Gamma(X;\mathcal{L}_X)$ (by hands) and called it the  
characteristic cycle of $F$. 
Later in \cite[Section 9.4]{KS2}, Kashiwara and Schapira redefined it functorially by
using their functor ``$\mu \mathrm{hom}$".
See \cite[Section 9.4]{KS2} for the details.
\begin{example}
For a submanifold $M\subset X$ and $\mathbf{k}_M \in \Dbrc(X)$
we have $\CCyc(\mathbf{k}_M) =[T_M^\ast X]$. 
\end{example}
Assume that the support of $F\in \Dbrc(X)$ is compact.
Then our objective here is to describe its global Euler-Poincar\'e index
\begin{equation*}
\chi (X; F) \coloneq \sum_{j\in \bbZ} (-1)^j \mathrm{dim}_{\mathbf{k}} H^j(X;F)
\end{equation*}
in terms of $\CCyc (F)$. 
For this purpose, for a continuous section $\sigma\colon X\longrightarrow T^\ast X$
of $\pi_X\colon T^\ast X\longrightarrow X$
let $[\sigma]\in H_{\sigma(X)}^0 (T^\ast X; \pi_X^!\mathbf{k}_X)$
be the element naturally defined by $\sigma$ (see \cite[Definition 9.3.5]{KS2}) and 
$\#\left( [\sigma]\cap \CCyc(F)  \right)$ its intersection number with 
the Lagrangian cycle $\CCyc(F)$ (see \cite[Section 9.3]{KS2} for the definition).
Then we have the following microlocal index theorem.

\begin{theorem}[Kashiwara \cite{Kas10} and Kashiwara-Schapira \cite{KS2}]\label{KSIND}
Assume that the support of $F\in \Dbrc(X)$, is compact.
Then for any continuous section $\sigma\colon X\longrightarrow T^\ast X$
of $\pi_X \colon T^\ast X\longrightarrow X$
we have
\begin{equation*}
\chi (X;F) =\# \left( [\sigma] \cap \CCyc(F) \right).
\end{equation*}
\end{theorem}

If $X$ is endowed with a Riemannian metric, we have a natural
isomorphism $T^\ast X\simeq TX$ and hence Theorem \ref{KSIND} for the constant sheaf
$\mathbf{k}_X\in \Dbrc(X)$ is nothing but the classical 
Poincar\'e-Hopf theorem for smooth vector fields on $X$.
Namely, Theorem \ref{KSIND} is a vast generalization of 
the Poincar\'e-Hopf theorem. 
In \cite[Section 9.4]{KS2} Kashiwara and Schapira proved that under some natural conditions
the characteristic cycles behave functorially with
respect to morphisms $f\colon Y\longrightarrow X$ of manifolds.
Later in \cite{SV1}, Schmid and
Vilonen relaxed their conditions. 
In \cite{Kas88} and \cite[Section 9.6]{KS2} Kashiwara and Schapira also obtained
Lefschetz fixed point formulas for $\mathbb{R}$-constructible sheaves. 
In \cite{K-Sch}, \cite{M-U-V}, \cite{SV2} and \cite{SV3} 
these results were effectively used to solve 
the conjectures proposed in 
Kashiwara's programs \cite{Kas88} and \cite{Kas93} 
of studying the admissible representations 
of real Lie groups via $\bbR$-constructible sheaves on 
flag manifolds (for some related results, see 
\cite{HK1}).  On the other hand, Nadler and Zaslow 
used characteristic cycles 
to construct the Fukaya categories of contangent bundles 
in \cite{Nad} and \cite{Na-Z}.  
For the further development of the theory of 
characteristic cycles and Lefschetz fixed 
point formulas, see e.g. \cite{B-M-M}, \cite{Gi}, \cite{GM2}, \cite{G-1}, \cite{G-3}, 
\cite{I-M-T}, \cite{KS4}, \cite{KT25}, \cite{M-T-2}, \cite{SS}, \cite{Schu}.

\section{Bernstein-Sato polynomials}
The theory of $\SD$-modules and constructible sheaves plays a prominent role 
also in the study of $b$-functions i.e. Bernstein-Sato polynomials and Milnor 
monodromies. First of all, let us recall the definition of $b$-functions. 
Let $X$ be a complex manifold and $f \colon X \longrightarrow \bbC$ a non-constant 
holomorphic function on it. Let $\SD_X[s]$ be the sheaf of (non-commutative) 
rings on $X$ of polynomials of $s$ with coefficients in $\SD_X$. 
Then it was proved by Bj\"{o}rk \cite{Bj2} that for a point 
$x_0 \in f^{-1}(0) \subset X$ there exists a ``non-zero" polynomial 
$b(s) \in \bbC [s]$ such that we have
\begin{equation}\label{funeq} 
P(s)f^{s+1} \equiv b(s)f^s \quad \mathrm{on} \ U \setminus f^{-1}(0)
\end{equation}
for an open neighborhood $U$ of $x_0$ in $X$ 
and some $P(s) \in \SD_X[s](U)$. We call the (non-zero) monic polynomial of 
the lowest degree satisfying this condition the (local) Bernstein-Sato 
polynomial or the (local) $b$-function of $f$ at $x_0 \in f^{-1}(0) \subset X$ 
and denote it by $b_{f,x_0}(s) \in \bbC [s]$.
Replacing $X$ by $U \subset X$ and suppressing the symbol $x_0 \in X$ we sometimes 
denote it simply by $b_f(s) \in \bbC [s]$. In the algebraic case i.e. if $f$ 
is a complex polynomial, the existence of $b$-functions was proved by Bernstein 
\cite{Ber2}. 
Bernstein-Sato polynomials are useful to study the analytic continuations of local 
zeta functions defined as follows. 
For a sufficiently small neighborhood $U \subset X$ of the point $x_0 \in f^{-1}(0) \subset X$, 
a test function $\varphi \in C^\infty_0 (U)$ i.e. a ($\bbC$-valued) 
$C^\infty$-function with compact support on it and $s \in \bbC$, the integral
\begin{equation*}
\zeta_f(\varphi)(s) \coloneq \int_{U}\vert f(x)\vert^{2s} \varphi (x)
dx_1 \wedge d \overline{x_1} \wedge 
\cdots \wedge dx_n \wedge d \overline{x_n}
\end{equation*}
converges if the real part $\Re s$ of $s$ is large enough. Using the above 
functional equation \eqref{funeq} in the definition of $b_{f,x_0}(s) \in \bbC[s]$ 
and integrations 
by parts, we can easily show that the function $\zeta_f(\varphi)(s)$ of $s$ can be 
continued to a 
meromorphic function on the whole complex plane $\bbC$ and its poles are contained 
in the set
\begin{equation*}
\bigcup_{i=0,1,2,\ldots}\{b_{f,x_0}^{-1}(0)-i\} \: \subset \bbQ \subset \bbC.
\end{equation*}

We call $\zeta_f(\varphi)(s)$ the local zeta function 
associated to $f$ and $\varphi$. 
Moreover in \cite{Kas4}, creating a new theory of direct images of $\SD$-modules, 
Kashiwara proved the following theorem.

\begin{theorem}[Kashiwara \cite{Kas4}]\label{kasbft}
All the roots of the $b$-function $b_{f,x_0}(s) \in \bbC[s]$ are negative rational numbers 
i.e. $b_{f,x_0}^{-1}(0) \subset \bbQ_{<0} \coloneq \{ \alpha \in \bbQ \,\vert\, \alpha <0\}$.
\end{theorem}

In order to explain the relation of the $b$-function $b_{f,x_0}(s)\in\bbC[s]$ with the 
Milnor monodromies of $f$, for a point $x\in f^{-1}(0)\subset X$ let 
$F_x\subset X\setminus f^{-1}(0)$ be the Milnor fiber of $f$ at $x$ and 
\begin{equation*}
\Phi_{j,x} \colon H^j(F_x;\bbC) \simto H^j(F_x;\bbC) \qquad (j\in \bbZ)
\end{equation*}
the Milnor monodromies of $f$ at $x$ (see \cite{Milnor} for the definitions). 
Then for the finite subset 
\begin{equation*}
E_{f,x_0} \coloneq \Biggl\{\lambda\in\bbC \,\Biggm\vert\,
\begin{aligned}
&\lambda \mathrm{\ is\ an\ eigenvalue\ of\ }\Phi_{j,x}\mathrm{\ for\ some\ }\\
&j\in\bbZ \mathrm{\ and\ } x\in f^{-1}(0) \mathrm{\ close\ to\ the\ point\ } x_0 
\end{aligned}
\Bigg\} \quad \subset\bbC
\end{equation*}
of $\bbC$ we have the following Kashiwara-Malgrange theorem. Recall that by the 
monodromy theorem of Landmann, Grothendieck et al. 
all the elements of $E_{f,x_0} \subset \bbC$ are roots of unity.

\begin{theorem}[Kashiwara \cite{Kas8} and Malgrange \cite{Ma2}, 
\cite{Ma3}]\label{thm:KashiwaraMalgrange}
We have an equality
\begin{equation*}
\{ \exp (2\pi \sqrt{-1} \alpha ) \ \vert \ \alpha \in b_{f,x_0}^{-1}(0)\}=E_{f,x_0}.
\end{equation*}
\end{theorem}
By Theorem \ref{thm:KashiwaraMalgrange}, for 
the test function $\varphi \in C^\infty_0 (U)$ if 
$\alpha\in\bbQ_{<0}$ is a pole of the local zeta function $\zeta_f(\varphi)(s)$ 
then we have $\exp(2\pi\sqrt{-1}\alpha) \in E_{f,x_0}$. 
Let $i_f \colon X \longhookrightarrow X \times \bbC_t \, \big(x\longmapsto(x,f(x))\big)$ 
be the graph embedding of $X$ into $X \times \bbC_t$ and 
$\SM\coloneq \mathrm{D} i_{f \ast}
\SO_X \in \Modrh (\SD_{X \times \bbC_t})$ the direct 
image of the regular holonomic $\SD_X$-module $\SO_X$ by $i_f$. Then Kashiwara and 
Malgrange proved Theorem \ref{thm:KashiwaraMalgrange} by introducing 
the so-called Kashiwara-Malgrange filtration 
$\{V_\alpha \SM\}_{\alpha\in\bbQ}$ of $\SM$ 
along $X\times \{0\} \subset X\times\bbC_t$
indexed by $\bbQ$ and defined as follows. 
Let $\SI \coloneq \SO_{X\times \bbC}t \subset \SO_{X\times \bbC}$ be the defining 
ideal of the complex hypersurface $X\simeq X\times\{0\}\subset X\times\bbC$ and for 
$j\in\bbZ$ set 
\begin{equation*}
V_j(\SD_{X\times\bbC})\coloneq \{P\in\SD_{X\times\bbC} \ \vert \ 
P\SI^k \subset \SI^{k-j} \ (k \in\bbZ) \} \quad \subset \SD_{X\times\bbC}.
\end{equation*}
Then locally we have 
\begin{equation*}
V_0(\SD_{X\times\bbC}) = \SO_{X\times\bbC}\langle \partial_{x_1},\ldots,\partial_{x_n},t\partial_t \rangle
\end{equation*}
and $V_{-1}(\SD_{X\times\bbC})=tV_0(\SD_{X\times\bbC})$, where $x=(x_1,\ldots,x_n)$ is a 
local coordinate of $X$. We thus obtain a filtration $\{V_j(\SD_{X\times\bbC})\}_{j\in\bbZ}$ 
of $\SD_{X\times\bbC}$ indexed by $\bbZ$. 
Let $\delta(t-f(x)) \in \SM= \mathrm{D} i_{f \ast} \SO_X$ be the canonical 
generator of $\SM$ and set $\Theta\coloneq -\partial_t t\in\SD_{X\times\bbC}$. 
Then for a 
polynomial $b(s)\in\bbC[s]$ the functional equation \eqref{funeq} 
is equivalent to the one 
\begin{equation*}
b(\Theta)\delta(t-f(x)) \in V_{-1}(\SD_{X\times\bbC})\delta(t-f(x)) 
\end{equation*}
on a neighborhood of $x_0 \in f^{-1}(0) \subset X \simeq 
X \times \{ 0 \} \subset X \times \bbC$ in $X \times \bbC$ 
(see \cite[Section 5]{CDM}). 
As $\SM$ is 
generated by $\delta(t-f(x)) \in \SM$ over $\SD_{X\times\bbC}$, 
by Theorem \ref{kasbft} for any section $m \in \SM$ of $\SM$ 
we can show that there exists a ``non-zero" polynomial 
$b(s) \in \bbQ [s]$ with coefficients in $\bbQ$ such that we have 
\begin{equation}
b( \Theta ) m \ \in V_{-1}(\SD_{X\times\bbC})m
\end{equation}
(see \cite{M-S}). 
We denote the (non-zero) monic polynomial of 
the lowest degree satisfying this condition by 
$b_{m}(s) \in \bbQ [s]$. Note that for 
the canonical generator $m_0 \coloneq \delta(t-f(x)) \in \SM$ of $\SM$ 
we have $b_{m_0}(s)=  b_{f,x_0}(s)$. 
Then for $\alpha \in \bbQ$ we set 
$\bbQ_{\leq \alpha} \coloneq \{ \beta \in \bbQ \ | \ 
\beta \leq \alpha \}$ and 
\begin{equation*}
V_\alpha \SM \coloneq \{ m \in \SM \ | \ 
b_m^{-1}(0) \subset \bbQ_{\leq \alpha} \}  
\qquad \subset \SM. 
\end{equation*}
It turns out that $V_\alpha \SM$ is a 
$V_{0}(\SD_{X\times\bbC})$-submodule of $\SM$. 
The Kashiwara-Malgrange filtration $\{V_\alpha \SM\}_{\alpha\in\bbQ}$ of $\SM$ 
thus obtained is related to the Milnor 
monodromies of $f$ as follows. 
First of all, recall that for Deligne's nearby cycle functor 
\begin{equation*}
\psi_f(\,\cdot\,) \colon \Db(X) \longrightarrow \Db\bigl(f^{-1}(0)\bigr)
\end{equation*}
introduced in \cite{De3} and a point $x\in f^{-1}(0)$ there exist isomorphisms 
\begin{equation*}
H^j\psi_f(\bbC_X)_x \simeq H^j(F_x;\bbC) \quad (j\in\bbZ)
\end{equation*}
(see e.g. \cite{Le-2} and \cite[Theorems 2.6 and 2.7]{Tak-2}). 
Let $\iota \colon f^{-1}(0) \longhookrightarrow 
X\simeq X\times\{0\}=t^{-1}(0)$ be the inclusion map. Then by the isomorphism 
\begin{equation*}
\DR_{X\times\bbC}(\SM) \simeq {i_f}_*(\bbC_X[\dim X])
\end{equation*}
and \cite[Proposition 4.2.11]{Di}, \cite[Exercise VIII.15]{KS2} 
we obtain an isomorphism 
\begin{equation*}
\psi_t\bigl(\DR_{X\times\bbC}(\SM)\bigr) \simeq \iota_*\psi_f(\bbC_X)[\dim X].
\end{equation*}
In \cite{Kas8} Kashiwara proved that for any $\alpha\in\bbQ$ the $\SD_X$-module 
$\gr_\alpha^V(\SM) \coloneq V_\alpha(\SM)/V_{<\alpha}(\SM)$ is regular holonomic and 
there exists an isomorphism 
\begin{equation*}
\psi_t\bigl(\DR_{X\times\bbC}(\SM)\bigr) [-1] \simeq 
\bigoplus_{-1 \leq \alpha <0} \DR_X(\gr_\alpha^V(\SM))
\end{equation*}
(see \cite{M-M} and \cite{M-S} for the details). For $\lambda\in\bbC^*$ 
if there exists $\alpha \in \bbQ$ such that $-1 \leq \alpha <0$ 
and $\lambda = \exp(2\pi\sqrt{-1}\alpha)$ we set 
\begin{equation*}
\psi_{t,\lambda} \bigl(\DR_{X\times\bbC}(\SM)\bigr) 
\coloneq  \DR_X(\gr_\alpha^V(\SM)) [1]. 
\end{equation*} 
Otherwise, we set it to be zero. Then we obtain a decomposition 
\begin{equation*}
\psi_t\bigl(\DR_{X\times\bbC}(\SM)\bigr) \simeq \bigoplus_{\lambda\in\bbC^*} 
\psi_{t,\lambda} \bigl(\DR_{X\times\bbC}(\SM)\bigr)
\end{equation*}
of $\psi_t\bigl(\DR_{X\times\bbC}(\SM)\bigr) \in \Db(X)$ with respect to 
the (generalized) eigenvalues $\lambda\in\bbC^*$ of its monodromy 
automorphism. 
Similarly, we have a decomposition
\begin{equation*}
\psi_f(\bbC_X) \simeq \bigoplus_{\lambda\in\bbC^*}\psi_{f,\lambda}(\bbC_X)
\end{equation*}
of $\psi_f(\bbC_X) \in \Db\bigl(f^{-1}(0)\bigr)$ and for the (generalized) eigenvalue 
$\lambda$-part $\psi_{f,\lambda}(\bbC_X) \subset \psi_f(\bbC_X) \, (\lambda\in\bbC^*)$ 
of the Milnor monodromy there exists an isomorphism 
\begin{equation*}
\psi_{t,\lambda}\bigl(\DR_{X\times\bbC}(\SM)\bigr) \simeq \iota_*\psi_{f,\lambda}(\bbC_X)[\dim X].
\end{equation*}
This implies that the left hand side is zero 
on a neighborhood of $x_0 \in f^{-1}(0) \subset X$ in $X$ 
if and only if $\lambda \notin E_{f,x_0}$. 
Moreover, when $\lambda\in\bbC^*$ is a root of unity i.e. 
there exists $\alpha \in \bbQ$ such that $-1 \leq \alpha <0$ 
and $\lambda = \exp(2\pi\sqrt{-1}\alpha)$, 
by the Riemann-Hilbert correspondence the left hand side is zero 
on a neighborhood of $x_0 \in f^{-1}(0) \subset X$ in $X$ 
if and only if $\gr_\alpha^V(\SM) \simeq 0$ there. 
As we always have $-1 \in b_{f,x_0}^{-1}(0)$ and 
$1 \in E_{f,x_0}$, for the proof of Theorem \ref{thm:KashiwaraMalgrange} 
it suffices to consider only the monodromy eigenvalues 
$\lambda \in E_{f,x_0}$ such that $\lambda \not= 1$. 
Then Theorem \ref{thm:KashiwaraMalgrange} follows from the equality 
\begin{equation*}
(b_{f,x_0}^{-1}(0) + \bbZ ) \setminus \bbZ = 
(b_{m_0}^{-1}(0) + \bbZ ) \setminus \bbZ = 
\{ \alpha \in \bbQ \setminus \bbZ \ | \ 
\gr_\alpha^V(\SM) \not= 0 \} 
\end{equation*}
obtained easily from the basic properties of 
the Kashiwara-Malgrange filtration $\{V_\alpha \SM\}_{\alpha\in\bbQ}$ of $\SM$. 
Morihiko Saito used 
the Kashiwara-Malgrange filtration to construct his theory of pure and mixed Hodge modules 
in \cite{Sa1} and \cite{Sa2}, generalizing that of 
the mixed Hodge structures of Deligne \cite{De2} 
(see \cite{EZ-2}, \cite{EZ-L}, \cite{P-S} and \cite{Schnell}). 
Hodge modules thus obtained have numerous applications in algebraic geometry and singularity 
theory. In particular, in \cite{Sa1} Saito proved a 
decomposition theorem for the direct images of pure polarizable 
Hodge modules by projective morphisms. Later in \cite{Kas98} 
Kashiwara conjectured that similar results hold true also 
for arbitrary semisimple holonomic $\SD$-modules. In the 
regular case, using also the results of Sabbah \cite{Sab3}, 
this conjecture was proved by Mochizuki in \cite{Mo-1} 
and \cite{Mo-2} (for related results, see also 
\cite{Dr} and \cite{Mo-4}). For the other applications of 
Hodge modules, see \cite{B-S}, \cite{D-L-2}, \cite{D-S-1}, \cite{D-S-2}, 
\cite{D-S-3}, \cite{Ke-Sc}, 
\cite{M-T-7}, \cite{M-O-P-W}, \cite{M-P}, \cite{Po-Sc1}, \cite{Po-Sc2}, 
\cite{Sab2}, \cite{Sab4}, \cite{Tan1} etc. 
For the further development of the theory of Bernstein-Sato polynomials, 
see \cite{A-G-L-B}, \cite{B-M-S}, \cite{B-V-W-Z}, 
\cite{E-L-S-V}, \cite{G}, \cite{K-book}, 
\cite{M-N}, \cite{Oa2}, \cite{Sabbah-1}, \cite{Sabbah-2}, \cite{Sa3}, 
\cite{S-K-K-O}, \cite{Tak-1}, \cite{U-C}, \cite{Wal} etc. 

To end this section, as a 
typical example of important problems in mathematics motivated by the Kashiwara-Malgrange 
theorem, we shall introduce the monodromy conjecture proposed by Denef and Loeser in \cite{D-L}. 
Assume that $f$ is a regular function on a smooth variety $X$ and let 
$\nu \colon Y\longrightarrow X$ be a proper morphism of smooth varieties inducing an 
isomorphism $Y\setminus\nu^{-1}\bigl(f^{-1}(0)\bigr) \simto 
X\setminus f^{-1}(0)$ and $\nu^{-1}\bigl(f^{-1}(0)\bigr) = (f\circ\nu)^{-1}(0) \,\subset Y$ 
is a simple normal crossing divisor in $Y$. Let $(f\circ\nu)^{-1}(0) = E_1\cup\cdots\cup E_p$ 
be the irreducible decomposition of $(f\circ\nu)^{-1}(0)$ and for each non-empty subset 
$I\subset \{1,2,\ldots,p\}$ set 
\begin{equation*}
E_I^\circ \coloneq \Bigl( \bigcap_{i\in I}E_i \Bigr) \setminus
\Bigl( \bigcup_{j \notin I}E_j \Bigr) 
\quad \subset(f\circ\nu)^{-1}(0).
\end{equation*}
For a local coordinate $x=(x_1,\ldots,x_n)$ of $X$ and $1\leq i\leq p$ let $m_i\geq 1$ 
(resp. $l_i\geq 0$) be the order of the function $f\circ\nu$ (resp. the differential form 
$\nu^*(dx_1\wedge\cdots\wedge dx_n))$ along the irreducible component $E_i\subset Y$ of
$(f\circ\nu)^{-1}(0)$. Then we have the following result of Denef and Loeser.

\begin{theorem and definition}[Denef-Loeser \cite{D-L}]
For a point $x_0\in f^{-1}(0) \subset X$ the rational function 
\begin{equation*}
Z_{f,x_0}(s) \coloneq \sum_{I \ne\phi} \chi\big(E_I^\circ \cap \nu^{-1}(x_0)\big)
\cdot\biggl\{\prod_{i\in I}\frac{1}{m_is+l_i+1}\biggr\} \qquad \in\bbQ (s)^*
\end{equation*}
of $s$ does not depend on the choice of the embedded resolution $\nu\colon Y\longrightarrow X$ 
of $f^{-1}(0)\subset X$. We call $Z_{f,x_0}(s)\in\bbC(s)^*$ the topological zeta function of $f$ 
at the point $x_0 \in f^{-1}(0) \subset X$.
\end{theorem and definition}

In \cite{D-L} Denef and Loeser proved this theorem by using Grothendieck's 
Lefschetz fixed point formula. 
Later, they developed a new theory of motivic zeta functions and reproved it more elegantly 
(see \cite{D-L-3}). 
Motivated by the works of Igusa on $p$-adic integrals and 
the Kashiwara-Malgrange theorem, then in \cite{D-L} 
Denef and Loeser proposed the following conjecture.
\begin{conjecture}[Monodromy conjecture, Denef-Loeser \cite{D-L}] 
\leavevmode
\begin{enumerate}
\item[\textup{(i)}] If $\alpha\in\bbQ_{<0}$ is a pole of $Z_{f,x_0}(s)$, 
then we have $b_{f,x_0}(\alpha)=0$,
\item[\textup{(ii)}] If $\alpha\in\bbQ_{<0}$ is a pole of $Z_{f,x_0}(s)$, 
then we have $\exp(2\pi\sqrt{-1}\alpha) \in E_{f,x_0}$.
\end{enumerate}
\end{conjecture}
Note that by the Kashiwara-Malgrange theorem (i) implies (ii). Despite a lot of contributions 
by many mathematicians, the monodromy conjecture was completely proved only in dimension 
two i.e. in the case $\dim X=2$. See \cite{Veys} for an excellent review on this fascinating 
and mysterious conjecture. Aiming at its proof, in \cite{D-L-2} and \cite{D-L-3} Denef and Loeser 
introduced the motivic Milnor fibers of $f$, which are the motivic reincarnations of the 
classical Milnor fibers. For this purpose, they used the primitive decompositions of 
nearby cycle sheaves obtained by Saito in \cite{Sa2} 
(see \cite[Section 6]{Tak-2} for the details). 
For the further development of 
motivic Milnor fibers and their applications to monodromies of complex polynomials, see 
e.g. \cite{E-T-1}, \cite{G-L-M}, \cite{M-T-7}, \cite{R-1}, \cite{R-2}, 
\cite{Saito}, \cite{Stapledon}, \cite{Tak-2}, \cite{T-T}.

\section{Irregular Riemann-Hilbert correspondence}
Recall that in the classical analysis of (complex) ordinary differential equations
on the complex plane $\mathbb{C}$ not only regular singular points 
but also irregular ones were studied precisely.
Unlike the case of regular singular points, 
around an irregular singular point of a linear ordinary differential equation
$Pu=0$ it is almost always impossible to describe its holomorphic solutions explicitly.
Instead, on each sufficiently narrow sector along it 
the asymptotic expansions of their holomorphic solutions can be calculated.
In the middle of the 19th century,
Stokes observed that the asymptotic expansion of a holomorphic solution may change
if we move the direction of the sector.
Namely, he discovered the so-called Stokes phenomenon.
In the middle of the 20th century, Hukuhara, Turrittin and others found methods
to calculate the formal structure of an ordinary differential equation 
at an irregular singular point and hence the asymptotic expansions 
of its holomorphic solutions.
We call it the Hukuhara-Turrittin theorem.
In order to describe the formal structure at an irregular singular point $p\in \mathbb{C}$,
we take a ramified covering (ramification) $\rho \colon \widetilde{U} \longtwoheadrightarrow U$
of an open neighborhood $U$ of $p$ in $\mathbb{C}$ and consider the pull-back of 
the differential equation by it (see e.g. \cite{Wasow} for the details). 

In the early 21th century, in \cite{Sab6} 
Sabbah proposed a program 
to obtain a higher-dimensional analogue of the Hukuhara-Turrittin theorem.
Let $X$ be a complex manifold and $\SM \in \Modh(\SD_X)$ a meromorphic connection on it 
along a divisor $D\subset X$.
Then in \cite{Sab6} Sabbah conjectured that there exists a proper 
morphism $f\colon Y\longrightarrow X$ of complex manifolds inducing an isomorphism $Y\setminus f^{-1}(D) \simto
X\setminus D$ such that the divisor $f^{-1}(D)\subset Y$ is normal crossing 
and for each point $p\in f^{-1}(D)$ of it there exists a ramification 
$\rho \colon \widetilde{U} \longtwoheadrightarrow U$ of its open neighborhood $U$ in $Y$
along the normal crossing divisor $f^{-1}(D) \cap U\subset U$ for which 
the meromorphic connection 
\begin{equation*}
\mathrm{D} \rho^\ast \left( (\mathrm{D} f^\ast \SM )|_U \right) \quad 
\in \Modh(\SD_{\widetilde{U}})
\end{equation*}
on $\widetilde{U}$ 
along the normal crossing divisor $\rho^{-1} (f^{-1} (D) \cap U)\simto f^{-1}(D)\cap U$
in $\widetilde{U}$ has a normal form (in the formal sense).
This conjecture was proved by Sabbah \cite{Sab6} in dimension 2 i.e.
in the case $\mathrm{dim} X=2$.
After some important contributions by Andr\'e in \cite{Andre}, 
it was finally settled by Kedlaya \cite{Ked0}, \cite{Ked} and Mochizuki \cite{Mo-3}.
This far-reaching result in mathematics can be considered as an analogue
for meromorphic connections of Hironaka's celebrated theorem on the resolutions of 
singularities of algebraic varieties. For its applications to 
the algebraic de Rham cohomology groups of integrable connections, 
see \cite{H-1}, \cite{H-R}, \cite{Hu-Te}, \cite{KT25}.  

On the other hand, also in the early 21th century, 
Kashiwara and Schapira introduced a new notion of ind-sheaves in \cite{KS2.5}  
to detect the growth orders of holomorphic solutions to holonomic $\SD$-modules.
For this purpose, they enlarged the category of usual sheaves to that of ind-sheaves as follows.
Let $M$ be a good topological space (whose topology is locally compact, Hausdorff, 
countable at infinity and of finite soft dimension)
and $\mathrm{Mod}(\bbC_M)$ the category of sheaves of $\bbC_M$-modules on it.
We denote by $\mathrm{Mod}_\mathrm{c}(\bbC_M) \subset \mathrm{Mod}(\bbC_M)$
the full-subcategory of $\mathrm{Mod}(\bbC_M)$ consisting of sheaves with compact
support.
Then, just like Schwartz's distributions are the elements of the topological dual
of the space of $C^{\infty}$-functions with compact support,
the ind-sheaves of Kashiwara and Schapira on $M$ are the objects of the ``dual category" 
of $\mathrm{Mod}_\mathrm{c}(\mathbb{C}_M)$ (i.e. the ind-objects of 
$\mathrm{Mod}_\mathrm{c}(\mathbb{C}_M)$) satisfying some conditions 
(see \cite[Chapter 1]{KS2.5} for the details). 
It turns out the category $\mathrm{I} \bbC_M$ of ind-sheaves on $M$ thus obtained 
is abelian and there exists a fully faithful embedding
\begin{equation*}
\iota_M \colon \mathrm{Mod}(\bbC_M) \longrightarrow \mathrm{I} \bbC_M
\end{equation*}
of abelian categories. 
In \cite{KS2.5} Kashiwara and Schapira also showed that if $M$ is a real analytic manifold 
$\mathrm{I} \bbC_M$ contains the 
full subcategory of subanalytic sheaves on $M$ i.e. 
the sheaves on the subanalytic site of $M$ 
\footnote{As was shown in \cite{Kas16}, such ind-sheaves 
are sufficient to treat the irregular Riemann-Hilbert correspondence.} 
and defined the ind-sheaf $\SD b_M^\rmt$ of tempered distributions on it. 
Then for a complex manifold $X$, by taking the Dolbeault complex of $\SD b_{X_\bbR}^\rmt$
on its underlying real analytic manifold $X_\mathbb{R}$
they defined an object $\SO_X^\mathrm{t}\in \Db(\mathrm{I} \bbC_X)$ in the derived category
$\Db(\mathrm{I} \mathbb{C}_X)$ of ind-sheaves on $X_\bbR$ to obtain a functor
\begin{equation*}
\Sol_X^\rmt(\,\cdot\,) \colon \Dbhol(\SD_X)^{\mathrm{op}} \longrightarrow \Db(\mathrm{I} \bbC_X)
\quad (\SM \longmapsto \rhom_{\SD_X}(\SM, \SO_X^\rmt)),
\end{equation*}
where $\Dbhol(\SD_X)\subset \Db(\SD_X)$ stands for the full subcategory of $\Db(\SD_X)$
consisting of objects whose cohomology sheaves are holonomic $\SD_X$-modules.
Unfortunately, this functor fails to be fully faithful 
\footnote{Indeed, by the results in \cite[Section 7]{KS03}, 
we see that for the two holonomic $\SD$-modules 
$\SM_1= \SD_{\bbC} \exp ( \frac{1}{x})$ and 
$\SM_2= \SD_{\bbC} \exp ( \frac{2}{x})$ on $\bbC_x$ such that 
$\SM_1 \not\simeq \SM_2$ we have 
$\Sol_{\bbC}^\rmt ( \SM_1) \simeq \Sol_{\bbC}^\rmt ( \SM_2)$.}, 
and does not allow us to obtain a generalization of the 
classical Riemann-Hilbert correspondence to 
(not necessarily regular) 
holonomic $\SD$-modules (for further developments of the study of 
ind-sheaves and the functor $\Sol_X^\mathrm{t}(\,\cdot\,)$, 
see \cite{Mor}, \cite{Pre2}, \cite{Pre1}). 
Nevertheless, by using the results of Kedlaya and Mochizuki and the idea of 
adding one variable in Tamarkin 
\cite{Tama08} (for its applications to symplectic geometry, 
see e.g. \cite{G-K-S}), 
in \cite{DK16} D'Agnolo and Kashiwara established such a Riemann-Hilbert correspondence 
for holonomic $\SD_X$-modules as follows. 
Let $\overline{\bbR} \coloneq \mathbb{R}\cup \{ +\infty,-\infty \} \simeq [-1, 1]$
be the natural compactification of $\mathbb{R}$.
Then, first they constructed a certain quotient category 
$\Eb(\rmI\bbC_X)$ of the derived category $\Db(\mathrm{I} \bbC_{X\times \overline{\bbR}})$
of ind-sheaves on $X\times \overline{\bbR}$.
An object of $\Eb(\mathrm{I} \mathbb{C}_X)$ is called an enhanced ind-sheaf on $X$.
Next they constructed a functor
\begin{equation*}
\Sol_X^{\mathrm{E}}(\,\cdot\,) \colon \Dbhol(\SD_X)^{\mathrm{op}} 
\longrightarrow \Eb(\mathrm{I} \mathbb{C}_X) \quad (\SM \longmapsto \rhom_{\SD_X}(\SM, \SO_X^{\mathrm{E}}))
\end{equation*}
by defining an object $\SO_X^{\mathrm{E}} \in \Eb(\mathrm{I} \mathbb{C}_X)$.
For $\SM\in \Dbhol(\SD_X)^{\mathrm{op}}$ we call $\Sol_X^{\mathrm{E}}(\SM) 
\in \Eb(\mathrm{I} \mathbb{C}_X)$ the enhanced solution complex of $\SM$.
They constructed also a functor
\begin{equation*}
\rhom^\mathrm{E} (\,\cdot\,, \SO_X^{\mathrm{E}}) \colon \Eb(\mathrm{I} \mathbb{C}_X)
\longrightarrow \Db (\SD_X)^{\mathrm{op}}.
\end{equation*}
Then the following decisive result was obtained in \cite{DK16}.

\begin{theorem}[Irregular Riemann-Hilbert correspondence, 
D'Agnolo-Kashiwara \cite{DK16}]
The functor $\Sol_X^{\mathrm{E}}(\,\cdot\,)$ is fully faithful.
Moreover, for any $\SM \in \Dbhol(\SD_X)^{\mathrm{op}}$ the natural morphism 
\begin{equation*}
\SM \longrightarrow \rhom^{\mathrm{E}}(\Sol_X^{\mathrm{E}}(\SM), \SO_X^{\mathrm{E}})
\end{equation*}
is an isomorphism.
\end{theorem}
Namely, the enhanced solution complex $\Sol_X^{\mathrm{E}}(\SM) 
\in \Eb(\mathrm{I} \mathbb{C}_X)$ of $\SM \in \Dbhol (\SD_X)^{\mathrm{op}}$
contains all the informations of $\SM$ such as its formal types and Stokes structures etc.
and we can reconstruct $\SM$ from it.
We call it the irregular Riemann-Hilbert correspondence. 
Note that in \cite{Sab5} Sabbah also established a 
Riemann-Hilbert correspondence for holonomic $\SD_X$-modules  
by using some Stokes filtrations of local systems as in the same way as the proof
of Shibuya-Malgrange theorem (inspired by Deligne) in the one dimensional case 
(see e.g. \cite{BV}). 
In the theory of \cite{DK16} and \cite{KS16}, we have a nice functoriality 
of the Riemann-Hilbert correspondence. 
For instance, let $f\colon Y \longrightarrow X$ be a morphism of complex manifolds.
Then for $\SM \in \Dbhol(\SD_X)^{\mathrm{op}}$ there exists an isomorphism
\begin{equation*}
\mathsf{E} f^{-1} \Sol_X^{\mathrm{E}}(\SM) \simto \Sol_Y^{\mathrm{E}}(\mathrm{D} f^\ast \SM),
\end{equation*}
where $\mathsf{E} f^{-1} \colon \Eb(\mathrm{I} \mathbb{C}_X) \longrightarrow 
\Eb(\mathrm{I} \mathbb{C}_Y)$
is a natural functor defined by $f$. 
For the further development of 
the theory of the irregular Riemann-Hilbert correspondences,
see e.g. \cite{DK19}, \cite{H-H}, \cite{Ito}, 
\cite{Kas16}, \cite{KS16}, \cite{Kuwa}, \cite{Mo22}, \cite{Te}. 

The irregular Riemann-Hilbert correspondence of \cite{DK16} 
have a lot of applications in the study
of Fourier transforms of holonomic $\SD$-modules defined as follows.
Let $X=\bbC_z^N$ be a complex vector space of dimension $N$ and $Y=\bbC_w^N$ 
its dual.
We regard them as algebraic varieties and use the notations $\SD_X$ and $\SD_Y$ 
for the sheaves of the rings of ``algebraic" differential operators on them.
Let $W_N\coloneq \bbC[z, \partial_z] \simeq \Gamma(X; \SD_X)$ 
and $W_N^\ast \coloneq \bbC[w, \partial_w] \simeq \Gamma(Y; \SD_Y)$
be the Weyl algebras on $X$ and $Y$ respectively.
Then there exists a ring isomorphism
\begin{equation*}
W_N \simto W_N^\ast \quad (z_i \longmapsto -\partial_{w_i}, \quad \partial_{z_i} \longmapsto w_i)
\end{equation*}
by which we can endow a left $W_N$-module $M$ with a structure of a left $W_N^\ast$-module.
We call it the Fourier transform of $M$ and denote it by $M^\wedge$. Recall also that for 
the affine algebraic varieties $X$ and $Y$ we have equivalences of categories 
\begin{equation*}
\Modcoh(\SD_X) \simeq \Mod_f(W_N), \quad \Modcoh(\SD_Y) \simeq \Mod_f(W_N^*)
\end{equation*}
obtained by taking global sections (see \cite[Propositions 1.4.4 and 1.4.13]{H-T-T}), 
where for a ring $R$ we denote by $\Mod_f(R)$ 
the category of finitely generated left $R$-modules. 
Thus, for a coherent 
$\SD_X$-module $\SM\in\Modcoh(\SD_X)$ we can define its Fourier transform 
$\SM^\wedge\in\Modcoh(\SD_Y)$. It turns out that there exists an equivalence of categories 
\begin{equation*}
(\,\cdot\,)^\wedge \colon \Modh(\SD_X) \simto \Modh(\SD_Y)
\end{equation*}
between the subcategories of (algebraic) holonomic $\SD$-modules 
(see \cite[Proposition 3.2.7]{H-T-T}). 
Although the definition of Fourier transforms of holonomic 
$\SD$-modules is so simple as such, 
it is very hard to know their basic properties,  
such as singular loci, holonomic ranks, exponential 
factors and Stokes structures etc. in general 
(for the classical results on the Fourier transforms of holonomic 
$\SD$-modules, see e.g. \cite{BE04}, \cite{Br}, \cite{BMV}, \cite{Dai00}, 
\cite{Gar04}, \cite{Ma4}, \cite{Ma5}, \cite{Mo10}, \cite{Mochi18}, 
\cite{Sab08}). Indeed, let 
\begin{equation*}
X\overset{p}{\longleftarrow}X\times 
Y\overset{q}{\longrightarrow}Y
\end{equation*}
be the projections. 
Then by a result of Katz-Laumon \cite{KL85}, 
for a holonomic $\SD_X$-module
$\SM \in \Modh (\SD_X)$ we have an isomorphism
\begin{equation*}
\SM^\wedge \simeq
\mathrm{D} q_\ast( \mathrm{D} p^\ast \SM 
\Dotimes \SO_{X\times Y}
e^{- \langle z, w \rangle }),
\end{equation*}
where $\mathrm{D} p^\ast, \mathrm{D} q_\ast, 
\Dotimes$ are
the operations for algebraic $\SD$-modules
and $\SO_{X\times Y}
e^{- \langle z, w \rangle }$ is the integral connection
of rank one on $X\times Y$ associated to the canonical
pairing $\langle \cdot ,  \cdot \rangle : X\times Y \longrightarrow \bbC$. 
As the connection $\SO_{X\times Y}
e^{- \langle z, w \rangle }$ is irregular, the regularity of 
$\SM \in \Modh (\SD_X)$ does not imply that of the Fourier transform 
$\SM^\wedge \in \Modh (\SD_Y)$. 
For the applications of the 
irregular Riemann-Hilbert correspondence to the study of Fourier transforms of 
holonomic $\SD$-modules, see e.g. 
\cite{B-H-H-S}, \cite{DHMS17}, \cite{DK17}, \cite{hoh22}, 
\cite{IT20a}, \cite{IT20b}, \cite{KS16-2}, \cite{KT23}, \cite{Tak-3} etc.

\end{document}